\documentclass{article}

\usepackage{amsmath,amsthm,amsfonts,amssymb}

\theoremstyle{plain}
\newtheorem{theorem}{Theorem}[section]
\newtheorem{lemma}[theorem]{Lemma}

\newtheorem{corollary}[theorem]{Corollary}
\newtheorem{proposition}[theorem]{Proposition}

\theoremstyle{definition}
\newtheorem{definition}[theorem]{Definition}

\theoremstyle{remark}

\newtheorem{remark}[theorem]{Remark}
\newtheorem{example}[theorem]{Example}

\DeclareMathOperator*{\colim}{colim}

\DeclareMathOperator*{\holim}{holim}

\usepackage{mathrsfs}
\usepackage{amsfonts}
\usepackage[all]{xy}
\usepackage{fancyhdr}
\usepackage[left=1in,right=1in,top=1in,bottom=1in]{geometry}

\DeclareMathOperator*{\Fiber}{Fiber}

\DeclareMathOperator*{\Tate}{Tate}
\DeclareMathOperator*{\cofiber}{cofiber}

\def \C              {\mathcal C}

\def \A              {\mathcal A}
\def \X              {\mathcal X}

\def \M              {\mathcal M}
\def \I              {\mathcal I}

\def \P              {\mathbb{P}}

\def \F              {{\mathcal F}_A}

\def \a              {\bf a}

\begin{document}

\title{On Triples, Operads, and Generalized Homogeneous Functors}
\author{Randy McCarthy and Vahagn Minasian}
\date{}

\maketitle

\begin{abstract}
\noindent
We study the splitting of the Goodwillie towers of functors in various settings. In particular,
we produce splitting criteria for functors $F: \A \to M_A$ from a pointed category with coproducts  
to $A$-modules in terms of differentials of $F$. Here $A$ is a commutative $S$-algebra. 
We specialize to the case when $\A$ is the category of $\a$-algebras for an operad $\a$ 
and $F$ is the forgetful functor, and derive milder splitting conditions in terms of 
the derivative of $F$. In addition, we describe how triples induce operads, and prove that, 
roughly speaking, a triple $T$ is naturally equivalent to the product of its Goodwillie layers
if and only if it is an algebra over its induced operad.

\vspace{10pt}

\noindent
MCS: 55P43, 18D50, 55P99
\end{abstract}

\section{Introduction}

One of the central results in Homological Algebra is 
the Hochschild-Kostant-Rosenberg (HKR) theorem (e.g. see Chapter 3 of 
~\cite{Loday}), which states that for a smooth algebra map $k \to A$, the Hochschild
homology coincides with the differential forms:
$$ HH_\ast (A) = \Omega^\ast _{A|k}.$$
This is a `splitting' result in its essence, as it is a consequence of the collapsing the 
Andr\'e-Quillen  fundamental spectral sequence for HH by constructing certain section maps. 

The splitting nature of the HKR theorem is more apparent in its topological analogue 
which we developed in~\cite{McMin}. There we work in the symmetric monoidal category of 
$S$-modules invented by Elmendorf et al in~\cite{EKMM}, and prove that under suitable 
conditions, the Topological Hochschild homology of a connective commutative $S$-algebra $A$
decomposes:
$$\P_A \Sigma TAQ(A) \stackrel{\simeq}{\to} THH(A),$$
where $\P _A (X) = \bigoplus X^{\wedge_A n}/\Sigma_n$ is the symmetric algebra triple, and 
$\Sigma TAQ(A)$ is the suspension of the Topological Andr\'e-Quillen homology.

The proof of this result is, in essence, a two step process. First, we develop splitting 
criteria for commutative $S$-algebras, then show that if the commutative $S$-algebra $A$
is sufficiently `nice' (or more precisely smooth in terminology of~\cite{McMin}), then the 
commutative $S$-algebra $THH(A)$ satisfies these splitting criteria. 
This served as one of the main motivations 
for the present work, as it gave rise to natural questions. Is the {\it commutative} algebra
structure special, or can the splitting results be extended to algebras of other types? 
More generally, when can we decompose functors? Addressing the ambiguity of the last question, 
we will explain shortly in what sense it is a generalization of the first one, but now 
we digress a little to make the notion of `decomposing functors' more precise, which leads 
to an explanation of the terminology used in the title.  

In a series of papers~\cite{Good1},~\cite{Good2} and ~\cite{Good3}, Tom Goodwillie developed a
theory (the calculus of homotopy functors), which has since become an important tool in 
homotopy theory. The central ingredient of his theory is a construction of a natural inverse 
limit system of functors
$$\cdots \to P_n F \to P_{n-1} F \to \cdots \to P_1 F$$
for a suitable functor $F$. This system, which is usually referred to as Goodwillie or Taylor tower 
of $F$, plays a role analogous to that of Taylor series in real variable calculus. 
In particular, often, though not always, the system $\{P_n F(X) \}$ converges to $F(X)$. 
However, unlike functions of real variables, the limit of the Taylor tower of functors does not 
necessarily decompose into the product of the fibers (or layers) 
$D_n F(X) \stackrel{def}{=} fiber [P_n F(X) \to P_{n-1} F(X)]$ (or equivalently, 
the Taylor tower does not necessarily split). For the cases when it does, we introduce a new 
terminology. 

To do so, first recall that by definition (due to T. Goodwillie), the functor $G$ is 
$n$-homogeneous if the Taylor polynomials $P_i G$ are contractible for all $i<n$ and 
$P_i G \simeq G$ for $i \geq n$. In particular, the layers $D_n F$ are $n$-homogeneous.
Thus, when the Goodwillie tower splits, its homotopy limit is equivalent to a product
of homogeneous functors. Keeping in mind that {\it generalized} Eilenberg-MacLane spaces
are precisely the spaces which decompose into products of fibers of their Postnikov 
systems, we present our definition.

\begin{definition}
Functor $F$ is called a {\it generalized} homogeneous functor if it is equivalent to a product 
of $n$-homogeneous functors (not necessarily for the same $n$). 
\end{definition} 

Functors with splitting Taylor towers provide the principal class of examples of generalized 
homogeneous functors, since for such functors, the functor 
$P_{\infty} F \stackrel{def}{=} \holim P_n F$ is equivalent to 
the product of the ($n$-homogeneous) fibers of the tower, and thus is generalized homogeneous.

We return to our motivational example of the HKR theorem developed in~\cite{McMin}. There, as a key 
step, we construct a tower of functors that approximates the forgetful functor $U$ from the 
category of non-unital commutative $A$-algebras to $A$-modules, and discuss conditions under 
which that tower splits. Here $A$ is a cofibrant commutative $S$-algebra and 
`approximates' means that under suitable conditions, the homotopy limit of the tower is 
equivalent to the functor. Not unexpectedly, this tower is equivalent to the Taylor tower 
of $U$, though we have intentionally avoided the Goodwillie Calculus language 
in~\cite{McMin} to increase its accessibility.

In this paper, we consider the conditions under 
which the forgetful functor from the category of algebras of other types (such as associative, 
or Lie) is generalized homogeneous at a given object. To make this statement precise we utilize 
the formalism 
provided by the language of operads. In other words, for a fixed operad $\a$, we consider the 
forgetful functor $U_{\a}: \C_{\a} \to \M_A$, where $\C_{\a}$ is the category of $\a$-algebras 
and $\M_A$ - of $A$ modules, and explore the question of the splitting of the Goodwillie tower 
of $U_{\a}$ evaluated at a fixed $\a$-algebra $X$.  

More generally still, we discuss the splitting of the Goodwillie tower of a functor 
$F: \A \to \M_A$ from {\it any} pointed category with coproducts. In fact, we produce conditions 
which are both necessary and sufficient for such towers to split. To explain these, we recall 
a generalization of the notion of derivative, which is introduced in Section 5 of~\cite{Randy}, 
and is also briefly described here, in Sections~\ref{sec:derdif} and~\ref{sect:splitfunc}.
First, the derivative itself is simply the functor  $P_1 F = D_1 F$, which is linear and 
comes equipped with a natural (derivative) map $F(X) \to D_1 F(X)$.
As with functions of real variables, one can define a notion of a directional derivative in 
Goodwillie Calculus to extend this. In other words, we can take the derivative of the functor 
$F$ at an object $Y$ in direction of an object $X$ (see Definition~\ref{def:diff}). 
The derivative $D_1 F(X)$ is simply the directional derivative at the base point $\ast$ 
in direction of $X$. 
  
An additional piece of notation will allow us to state one of our main results in a rather 
compact form. Define $\X$ to be the full subcategory of $\A$ whose objects are $\vee _n X$ for 
all $n \geq 0$ (with $\vee_0 X =\ast$), and denote by $F|_{\X}$ the restriction of $F$ to 
$\X$.

\begin{theorem}
\label{theor:int1}
The functor $P_{\infty} (F|_{\X})$ is generalized homogeneous if and only if the natural 
(derivative) map from $F(X \vee Y)$ to the directional derivative of $F|_{\X}$ at $Y$ in 
direction of $X$ has a section for all $Y \in \X$. 
\end{theorem} 

In particular, this theorem provides criteria under which the forgetful functor 
$U_{\a}$ evaluated at an algebra $X$ (and hence the algebra $X$ itself) decomposes 
(if the homotopy limit of the Taylor tower of $U_{\a}$ is equivalent to $U_{\a}$ itself).  
This approach however ignores the 
additional structure of an algebra that the objects under discussion possess. 
The special case of commutative algebras mentioned above though, suggests that it is 
of critical importance. Indeed, in~\cite{McMin}, we showed that the Goodwillie tower
of the forgetful functor $U$ from the category of commutative $A$-algebras to $A$-modules
splits if the natural derivative map $U(X) \to D_1 U(X)$ has a section. In other words, 
to produce a decomposition for commutative algebras, we required that 
in terminology of Theorem~\ref{theor:int1}, the directional derivative be 
equipped with a section only for $Y = \ast$, as opposed to $Y=\vee_n X$ for all $n$.

We are able to get a similar result for algebras over other operads as well. We do restrict 
however to operads $\a$ with ${\a}(1)$ 
equivalent to the unit $A$ of the symmetric monoidal category $\M_A$. 
Note that this is not a 
very restrictive assumption as most operads naturally occurring in literature satisfy it. 

\begin{theorem}
\label{theor:int2}
Let $\a$ be an operad with ${\a}(1)$ equivalent to $A$, and let $C$ be an $\a$-algebra such 
that the natural 
derivative map $U_{\a} (C) \to D_1 U_{\a}(C)$ has a section in the category of $A$-modules. 
Then  
$$\holim_n P_n U_{\a} (C) \simeq \prod {\a}(n) \wedge_{\Sigma_n} [D_1 U_{\a} (C)]^{\wedge n}.$$ 
\end{theorem}

Note that as a trade off for relaxing the splitting criteria of 
Theorem~\ref{theor:int1}, we claim that the Goodwillie tower of $U_{\a}$ splits only at $C$ 
and not all finite multiples of $C$. 

The special case of the commutative algebra operad ${\bf e_{\infty}}$ implies 
the Hochschild-Kostant-Rosenberg theorem. Another immediate application of this result recovers
the theorem of Leray on the structure of commutative quasi Hopf algebras, and consequently the
Poincar\'e-Birkhoff-Witt theorem. 

Yet another special case of Theorem~\ref{theor:int1} is obtained by considering only 
those functors which have the additional structure of a triple. In fact, a splitting 
result on a triple $T$ would be a natural culmination of the two main theorems presented so 
far, since for all $X$, $T(X)$ is equipped with a structure of an algebra. Complications
arise due to the fact that not every triple is induced by an operad, and thus $T(X)$ may 
not be an algebra over a triple.
Hence Theorem~\ref{theor:int2} may not be applicable. 

To remedy this problem, we present a construction which is of interest on its own. It is 
based on a simple observation that if $T_{\a}$ is the triple associated with the operad 
$\a$ in the symmetric monoidal category $\M_S$ of $S$-modules, then we can recover the 
$n$'th space ${\a}(n)$ of the operad by multilinearizing the functor 
$cr_n T_{\a}$ at each of its $n$ variable, where $cr_n$ is the $n$'th cross effect 
(see Section~\ref{sec:derdif} for a definition). In other words,
$${\a}(n) \simeq D_1 ^{(n)} cr_n T_{\a} (S, \cdots , S),$$
where $D_1 ^{(n)}$ indicates that we have applied $D_1$ successively with respect to each 
of the $n$ variables, and $S$ is the sphere spectrum, which is the unit of our symmetric 
monoidal category $M_S$.   

It turns out that this construction produces an operad even if we replace the triple 
$T_{\a}$ by a triple which is not necessarily associated with an operad. More precisely,
we show that for every triple $T$, we can define an operadic multiplication on the 
sequence of objects $\{ {\a}_T (n) \}$ given by 
${\a}_T (n) \stackrel{def}{=} D_1 ^{(n)} cr_n T (S, \cdots , S)$. We refer to this operad
as the operad induced by the triple $T$. Of course, it is not necessarily the case that the 
triple $T_{{\a}_T}$ associated with the operad ${\a}_T$ is equivalent to $T$. 
For all $X$, $T(X)$ is always a $T$-algebra, and,  when  $T$ and $T_{{\a}_T}$ are different, 
$T(X)$ may in addition have a structure of an  ${\a}_T$-algebra,
thus equipping $T(X)$ with two {\it different}
algebraic structures - arising from the triple $T$ and the operad ${\a}_T$. 

Exploration of these structures is at the root of understanding the splitting of the 
Goodwillie tower of the triple $T$, as we prove the following theorem.

\begin{theorem}
The Goodwillie tower of a triple $T$ in $\M_A$ splits at $X$ if and only if $T(X)$ is an 
${\a}_T$-ring spectrum, and the two algebra structures are compatible in some natural 
sense. 
\end{theorem}

The three splitting results described above form the core of this work.  
 
The paper is structured as follows. In Section~\ref{sec:derdif}, we briefly recall some basic 
definitions 
from Goodwillie Calculus, as well as set up the notation. Theorem~\ref{theor:int1} is 
proved in Section~\ref{sect:splitfunc}. Section~\ref{sect:algbar} is devoted to 
adopting P.May's two sided bar construction 
to the category $\C_{\a}$ of $\a$-algebras. In addition, geometric realizations and 
closed model structures on $\C_{\a}$ are discussed. In Section~\ref{sect:forget}, 
we introduce the forgetful 
functor $U_{\a}:\C_{\a} \to \M_A$ and compute its layers $D_n$ in terms of $D_1$. 
Also, for a special class of operads (to which we refer as `primitively generated operads'), 
we give an algebraic description of the Taylor 
tower of $U_{\a}$, which does not assume familiarity with Goodwillie Calculus. 
We employ these computations to prove Theorem~\ref{theor:int2} (as well as its 
`Calculus free' analogue for the special class of primitively generated operads) in 
Section~\ref{sect:splitalg}. There, we also discuss how the classical theorems of 
Leray and Hochschild,Kostant and Rosenberg can be recovered from our results.  
In Section~\ref{sect:tripop}, 
we show how a triple $T$ induces an operad and discuss some examples. Section~\ref{sect:splittrip} 
explores the two algebra structures on $T(X)$ (produced by the triple itself and the induced 
operad), and uses these to derive necessary and sufficient conditions for $T$ to be generalized 
homogeneous. Finally, in Section~\ref{sect:last}, we prove the technical results presented 
(without proof) in Section~\ref{sect:forget}.   

{\bf Acknowledgments.}
We would like to thanks Paul Goerss for introducing us to a number of papers that helped us 
to sort through model theoretic questions on the categories of algebras over operads.

\section{Derivatives and Differentials}
\label{sec:derdif}
Let $A$ be a cofibrant commutative $S$-algebra, and let $\M_A$ be the category of 
$A$-modules. In this section we give a brief summary of relevant (for our purposes) 
constructions of Goodwillie Calculus for the category of functors from a 
pointed category $\C$ to $\M_A$. It is no more than a restatement of the appropriate
constructions from~\cite{Randy}, which in turn
is an application of ~\cite{Good1}, ~\cite{Good2} and ~\cite{Good3}. 

We begin by recalling the definition of the cross effects of a functor 
$F: \C \to \mathcal{A}$, where $\C$ is a basepointed category with finite coproducts 
and $\mathcal{A}$ is an Abelian category. We denote the basepoint of $\C$ by $0$.
The notion of cross effects dates back to S.Eilenberg and S.MacLane (~\cite{EilMac}).
\begin{definition}
The n-th cross effect of $F$ is the functor $cr_n F: \C ^{\times n} \to \mathcal{A}$ 
defined inductively by 

\noindent
$cr_1 F(M) \oplus F( \ast) \cong F(M)$

\noindent
$cr_2 F(M_1, M_2) \oplus cr_1 F(M_1) \oplus cr_1 F(M_2) \cong cr_1 F(M_1 \vee M_2)$

\noindent
and in general,

\noindent
$cr_n F(M_1, \cdots M_n) \oplus cr_{n-1} F(M_1,M_3 \cdots M_n)  
\oplus cr_{n-1} F(M_2,M_3 \cdots M_n)$

\noindent
is equivalent to 
$$cr_{n-1} F(M_1 \vee M_2,M_3 \cdots M_n).$$
\end{definition}

To ease the notation we will denote the $n$-multifunctor $cr_n F(M, \cdots, M)$ by $cr_n F(M)$.

\begin{definition}
Given a functor $F$ from $\C$ to $\M_A$, we say that $F$ is degree $n$ if $cr_{n+1} F$
is acyclic as a functor from $\C^{\times n+1}$ to $\M_A$. That is, $cr_{n+1} F$ is 
contractible when evaluated on a collection of $n+1$ objects in $\C$. 
\end{definition}

Denote the category of functors of $n+1$ variables from $\C$ to $\mathcal{A}$ that
are reduced in each variable by $Func_\ast (\C^{\times n+1},\mathcal{A})$. Let 
$\bigtriangleup ^\ast$ be the functor from  $Func_\ast (\C^{\times n+1},\mathcal{A})$
to $Func_\ast (\C, \mathcal{A})$ obtained by composing a functor with the diagonal
functor from $\C$ to $\C^{\times n+1}$. When $\mathcal{A}$ is the category of chain
complexes over a commutative ring, the $(n+1)$st cross effect is the right adjoint
to $\bigtriangleup ^\ast$. Consequently, the functor  
$\perp _{n+1} \stackrel{def}{=}  \bigtriangleup ^\ast \circ cr_{n+1}$ is a cotriple on 
$Func_\ast (\C, {\mathcal A})$. For the case ${\mathcal A} = \M_A$,  $\bigtriangleup ^\ast$ 
and $cr_{n+1}$ are not adjoint, however the composite functor 
$\perp _{n+1} = \bigtriangleup ^\ast \circ cr_{n+1}$ is still a cotriple on $Func_\ast (\C, \M_A)$
(see Appendix of~\cite{Mauer} for a proof). 

\begin{definition}
Let $\perp _{n+1} =  \bigtriangleup ^\ast \circ cr_{n+1}$ be the above cotriple on 
$Func_\ast (\C, \M_A)$. For each $F \in Func_\ast (\C, \M_A)$, denote by 
$\perp_{n+1} ^{\ast +1} F$ the simplicial object whose simplices are $\perp_{n+1} ^{\ast +1}F$.
We define $P_n$ to be the functor from 
$Func_\ast (\C, \M_A)$ to $Func_\ast (\C, \M_A)$ given by
$$P_n F(X) = hocofiber [|\perp_{n+1} ^{\ast +1} F(X)| \stackrel{\epsilon}{\to} id(F(X))]$$
where $|\perp_{n+1} ^{\ast +1} F(X)| $ is the geometric realization of the simplicial object. 
Furthermore, let $p_n:id \to P_n$ be the natural transformation obtained from the 
homotopy cofiber.
\end{definition}
We note that the functor $P_n F(X)$ is degree $n$. In addition, if $F$ is already 
of degree $n$, then $p_n : F \to P_n$ is an equivalence.

Next we produce a natural transformation $q_{n}:P_{n} \to P_{n-1}$. 

Observe that we have the following formula relating the $n$'th and $n+1$'st 
cross effects:
\begin{equation}
\label{eq:cross}
cr_{n+1} F(X_1, \cdots, X_{n+1}) = 
cr_2 (cr_n F(X_1, \cdots, X_{n-1},-))(X_n,X_{n+1}).
\end{equation}
Fix an object $X$ of $\C$. Let $G(Y)= cr_n F(X, \cdots, X,Y)$. Then using the fold
map, we have $cr_2 G(X,X) \to G(X \vee X) \to G(X)$, which in turn gives us
$$cr_2 (cr_n F(X, \cdots, X,-))(X,X) \to cr_n F(X, \cdots, X,X).$$
Combine this map with the Equation~\ref{eq:cross} to produce a map
$ cr_{n+1} F(X, \cdots, X) \to cr_{n} F(X, \cdots, X)$, which induces the 
desired map $q_n$. 
 
\begin{definition}
The $n$'th layer or the $n$'th derivative of $F$ is the functor
$$D_n F(-) (X) \stackrel{def}{=} hofiber(q_n)(X).$$
In particular, the functor $D_1 F$ (which, by definition, is equivalent to $P_1 F$) will 
be referred to as ``the derivative'' or ``the linearization'' of $F$.  
\end{definition}
The linearization functor $D_1$ plays a rather central role, since the higher derivatives 
can be expressed in terms of $D_1$. So following~\cite{Randy} we introduce the following 
notation. For a 
functor $F$ of $n$ variables, we denote by $D_1 ^i F(X_1, \cdots , X_n)$ the derivative of $F$ 
obtained by holding all but $i'th$ variable constant, in other words $D_1$ is applied to 
the single variable functor obtained by fixing all but $i'th$ variables. In addition, denote 
by $D_1 ^{(n)}F$ the multilinearized functor $D_1 ^n \cdots D_1 ^2 D_1 ^1 F$.

The following proposition, which provides a description of objects in Goodwillie towers in 
terms of $D_1$, is the combination of Proposition 3.9 and a special case of Proposition 3.1
of~\cite{Randy}.

\begin{proposition}
\label{prop:deriv}
Let $F:\C \to \M_A$ be a functor into the category of $A$-modules. Then

{\bf 1.} $D_n F$ is naturally equivalent to $D_1 ^{(n)} cr_n F_{h\Sigma_n}$,   

{\bf 2.} $P_n cr_n F(X)$ is naturally equivalent to $D_1 ^{(n)} cr_n F(X)$.
\end{proposition}

We conclude this section with the definition of a {\it differential} of a functor, which mimics 
the notion of the directional derivative of a differentiable multi-variable function. 
See Section 5 of~\cite{Randy} for details.

\begin{definition}
\label{def:diff}
Let $F:\C \to \M_A$ be a functor, and $X$ and $Y$ be objects in $\C$. The differential 
of $F$ is the bifunctor defined by
$$
\nabla F(X; Y) = D_1 ^X [\Fiber (F(Y \vee X) \to F(Y))],
$$
where the superscript in $D_1 ^X$ indicates that the derivative is taken with respect to the 
variable $X$. Consequently, $\nabla F(X; Y)$ is linear in $X$ but not necessarily
in $Y$.
\end{definition}
 
As a simple example, observe that the differential at $Y=0$ is 
$$
\nabla F(X; 0)= D_1 ^X [\Fiber (F(0 \vee X) \to F(0))] \cong D_1 ^X cr_1 F(X) = D_1 F(X).
$$

\section{Splitting of Goodwillie Towers}
\label{sect:splitfunc}

For a functor $F$ into the category $\M _A$ of $A$-modules, the data of the previous 
Section~\ref{sec:derdif} assembles into a diagram (commonly referred to as the Goodwillie 
tower of $F$ at $X$)
$$
\xymatrix{
&
\cdots \ar[d]^{q_{n+1}}
& \\
F(X) \ar[r]^{p_n} \ar[rd]^{p_{n-1}}
&
P_n F(X) \ar[d]^{q_{n}} 
&
D_n F(X) \ar[l]_{d_n}\\
&
P_{n-1} F(X) \ar[d] 
&
D_{n-1} F(X) \ar[l]_{d_{n-1}}\\
&
\cdots
&
}$$ 
such that the homotopy limit $P_{\infty} F(X)$ of the inverse limit system 
$\{ P_n F(X) \}$ is equivalent to $F(X)$ under suitable conditions. While the 
convergence issue will get addressed at various points throughout this work, here
we concentrate on the question of developing conditions under which the above tower
splits. 

In this section alone, we will assume that the above tower does converge to $F(X)$. 
We do so to ease the notation while discussing the splitting problems. Naturally,
once the tower is constructed, it will split or not split regardless of the object 
to which it converges. In fact, in what follows we could simply replace $F(X)$ by 
$P_{\infty} F(X)$ to eliminate the issue of convergence.

We begin by making a few simple observations that will lead to necessary conditions
for splitting of the tower and will also simplify our terminology and arguments at later stages.  
Note that if the tower splits at $X$ (i.e. $F(X) \simeq \Pi D_n F(X)$), the derivative 
map $F(X) \to D_1 F(X)$ has a section. More generally, by Proposition~\ref{prop:deriv}
we have that 

\begin{equation}
\label{eq:eq1}
D_n F(X) \simeq D_1 ^{(n)}cr_n F(X)_{h\Sigma_n} \cong D_1 ^{(n)} F(\vee _n X)_{h\Sigma_n}.
\end{equation}
To see the second equivalence, recall that 
$$cr_n F(M_1, \cdots , M_n) \vee cr_{n-1} F(M_1,M_3, \cdots ,M_n)  
\vee cr_{n-1} F(M_2,M_3, \cdots ,M_n)$$
is isomorphic to 
$cr_{n-1} F(M_1 \vee M_2,M_3, \cdots ,M_n)$. However note that the terms 
$cr_{n-1} F(M_1,M_3, \cdots ,M_n)$
and $cr_{n-1} F(M_2,M_3, \cdots ,M_n)$ are constant in variables $M_2$ and $M_1$ 
respectively. Hence they vanish after we apply $D_1 ^2$ and $D_1 ^1$. Consequently,
$D_1 ^{(n)}cr_n F(M_1, \cdots , M_n) \simeq D_1 ^{(n)} 
cr_{n-1} F(M_1 \vee M_2,M_3, \cdots ,M_n)$. Iterating this construction, we get the second 
part of Equation~\ref{eq:eq1}.  
Further observe that 
\begin{equation*}
D_1 ^1 F(\vee _n X) = D_1 ^1 F(X \vee \vee_{n-1} X) \simeq 
D_1 ^1 [F(X \vee \vee_{n-1} X) - F(\vee_{n-1} X)] = \nabla F(X; \vee_{n-1} X),
\end{equation*}
where by $F(X \vee \vee_{n-1} X) - F(\vee_{n-1} X)$ we mean the fiber of the obvious map 
$F(X \vee \vee_{n-1} X) \to F(\vee_{n-1} X)$ which has a section. We will use 
this notation in the future. (See Definition~\ref{def:diff}.) In particular, the question 
of the existence of a splitting of the derivative map $F(X) \to D_1 F(X)$ is equivalent 
to the existence of a splitting of the map
$F(X) \to \nabla F(X; 0)$. 

Our immediate objective is to show that, 
more generally, if the Goodwillie tower of $F$ splits at $\vee_k X$ for all $k \geq 1$, then
the natural maps 
\begin{equation}
\label{eq:added1}
F(\vee_k X) = F(X \vee \vee_{k-1} X) \to \nabla F(X; \vee_{k-1} X) 
\end{equation}
also split. 
In other words, assuming that the Goodwillie tower of $F$ at $\vee_k X$ decomposes into a 
product, we need to construct splittings
\begin{equation*}
\nabla F(X; \vee_{k-1} X) \simeq \nabla \Pi D_n F(X; \vee_{k-1} X) 
\to \Pi D_n F(\vee_k X) \simeq F(\vee_k X). 
\end{equation*} 
To do so we analyze the individual factors on the two sides. Observe that 
\begin{eqnarray}
\label{eq:dnfa}
\lefteqn{
D_n F(\vee_k X) \simeq D_1 ^{(n)} cr_n F(\vee_k X)_{h\Sigma_n} =
D_1 ^{(n)} cr_n F(\vee_k X,\vee_k X, \cdots, \vee_k X )_{h\Sigma_n}}\\
& \simeq &
[D_1 ^{(n)} cr_n F(S,\vee_k X, \cdots, \vee_k X )\wedge (\vee_k X)]_{h\Sigma_n} 
\simeq \cdots \simeq D_1 ^{(n)} cr_n F(S) \wedge_{h\Sigma_n}(\vee_k X)^{\wedge n},\nonumber
\end{eqnarray}  
where $S$ is the sphere spectrum, and the last $n$ equivalences are by linearity of 
the derivative $D_1$. See ~\cite{Randy2} for more on this. On the other hand,
\begin{eqnarray*}
\lefteqn{
\nabla D_n F(X; \vee_{k-1} X) = 
D_1 ^1 [D_n F(X \vee \vee_{k-1} X)- D_n F(\vee_{k-1} X)]}\\
& \simeq &
D_1 ^1 [D_1 ^{(n)} cr_n F(S) \wedge_{h\Sigma_n} (\vee_k X)^{\wedge n} -
D_1 ^{(n)} cr_n F(S) \wedge_{h\Sigma_n} (\vee_{k-1} X)^{\wedge n}]\\ 
& \simeq &
D_1 ^1 [D_1 ^{(n)} cr_n F(S) \wedge_{h\Sigma_n} ((\vee_k X)^{\wedge n} 
-(\vee_{k-1} X)^{\wedge n})]\\
& \simeq &
D_1 ^1 [D_1 ^{(n)} cr_n F(S) \wedge_{h\Sigma_n} (n X \wedge (\vee_{k-1} X)^{\wedge n-1}
\vee \frac{n(n-1)}{2} X^{\wedge 2} \wedge (\vee_{k-1} X)^{\wedge n-2} \vee \cdots \vee
X^{\wedge n})]\\
& \simeq &
D_1 ^1 [D_1 ^{(n)} cr_n F(S) \wedge_{h\Sigma_n} (n X \wedge (\vee_{k-1} X)^{\wedge n-1})].
\end{eqnarray*}
To obtain the last equivalence we made use of the fact that for $l \geq 2$, 
$D_1 ^1 (X^{\wedge l} \wedge (\vee_{k-1} X)^{\wedge n-l})$ is contractible by 
Proposition 3.1 of~\cite{Randy} since $X^{\wedge l} \wedge (\vee_{k-1} X)^{\wedge n-l}$ is 
$l$-multireduced in the variable with respect to which the derivative is taken. 
Moreover, note that
\begin{equation*}
D_1 ^1 [D_1 ^{(n)} cr_n F(S) \wedge_{h\Sigma_n} (n X \wedge (\vee_{k-1} X)^{\wedge n-1})]
\simeq
D_1 ^{(n)} cr_n F(S) \wedge_{h\Sigma_n} [\vee_n X \wedge (\vee_{k-1} X)^{\wedge n-1}],
\end{equation*} 
to conclude that 
$\nabla D_n F(X; \vee_{k-1} X) \simeq 
D_1 ^{(n)} cr_n F(S) \wedge_{h\Sigma_n} [\vee_n X \wedge (\vee_{k-1} X)^{\wedge n-1}]$. 
We can restate the Map~\ref{eq:added1} in terms of this description and the description 
provided by Equation~\ref{eq:dnfa}. It is simply the map induced by projections
$$(\vee_k X)^{\wedge n} \cong (X \vee \vee_{k-1} X)^{\wedge n} \cong 
X^{\wedge n} \vee (\vee_n X \wedge (\vee_{k-1} X)^{\wedge n-1}) \vee \cdots 
\to \vee_n X \wedge (\vee_{k-1} X)^{\wedge n-1}.$$
Thus the desired splitting 
\begin{equation*}
\nabla D_n F(X; \vee_{k-1} X) \to D_n F(\vee_k X) 
\end{equation*}
is produced by inclusions $\vee_n X \wedge (\vee_{k-1} X)^{\wedge n-1} \to (\vee_k X)^{\wedge n}$. 
Further, observe that the projections $\Pi D_n F(Y) \to D_n F(Y)$ give rise to 
morphisms $\nabla \Pi D_n F(X;\vee_{k-1} X) \to \nabla D_n F(X;\vee_{k-1} X)$, 
which in turn induce a map into the product 
\begin{equation*}
\nabla \Pi D_n F(X;\vee_{k-1} X) \to \Pi \nabla D_n F(X;\vee_{k-1} X),
\end{equation*}
producing a composite
\begin{equation*}
\nabla F(X;\vee_{k-1} X) \simeq \nabla \Pi D_n F(X;\vee_{k-1} X) \to
\Pi \nabla D_n F(X;\vee_{k-1} X) \to \Pi D_n F(\vee_k X) \simeq F(\vee_k X),
\end{equation*}
which is the desired splitting.

In fact, the existence of these splittings is also a sufficient 
condition for the towers at $\vee_k X$ for all $k \geq 1$
to split. More precisely, we prove the following theorem. 

\begin{theorem}
\label{theor:splitfunc}
Let $F:\C \to \M_A$ be a functor from a pointed category $\C$ to the category of 
$A$-modules, and let $X$ be an object in $\C$. Then the following two statements 
are equivalent. 

{\bf 1.} $F$ and $X$ are such that for all $n$, in addition to the natural derivative maps 
$F( \vee _n X) \to \nabla F(X; \vee_{n-1} X) \cong D_1 ^i F( \vee _n X)$ 
(for $1 \leq i \leq n$), there are $n$ morphisms (one for each copy of $X$ in $\vee _n X$) 
$\nabla F(X; \vee_{n-1} X) \to F( \vee _n X)$, such that the $n$ composites 
\begin{equation}
\label{eq:split}
\nabla F(X; \vee_{n-1} X) \to F( \vee _n X) \to 
\nabla F(X; \vee_{n-1} X)
\end{equation} 
are equivalences. 

{\bf 2.} The Goodwillie tower of $F$ at $\vee_n X$ splits for all $n \geq 1$.   
\end{theorem}

\begin{proof}
Given our discussion preceding the theorem, we only need to show that the second statement 
follows from the first one. 

First we show that under conditions of Statement {\bf 1}, the Goodwillie tower of $F$ at $X$ 
splits.

To begin, observe that for $n=1$ the Equation~\ref{eq:split} simply states that the derivative 
map $F(X) \to D_1 F(X)$ has a section  $s:D_1 F(X) \to F(X)$, since 
$\nabla F(X; 0) \cong D_1 F(X)$. Composing the section $s$ with 
$p_2:F(X) \to P_2 F(X)$, we get a morphism $D_1 F(X) \to P_2 F(X)$. 

More generally, we would like to produce similar maps for higher degrees, in other words, we
are looking to construct sections $D_n F(X) \to P_{n+1} F(X)$ to $q_n$ for all $n >1$. To do so, 
consider the map 
$$D_1 ^1 F(\vee _n X) \simeq \nabla F(X; \vee_{n-1} X) \to F(\vee _n X)$$
and apply $D_1 ^2$ to the two sides to get
$D_1 ^2 D_1 ^1 F(\vee _n X)\to D_1 ^2 F(\vee _n X)$. Recalling that by 
assumption we have a map $ D_1 ^2 F(\vee _n X) \cong \nabla F(X; \vee_{n-1} X) \to F(\vee _n X)$, 
we produce a composite
$ D_1 ^2 D_1 ^1 F(\vee _n X)\to F(\vee _n X)$. Iterating this construction, i.e. applying 
$D_1 ^3, \cdots , D_1 ^n$ in succession, and composing the resulting morphism with the 
`$+$'-map, we get 
\begin{equation}
\label{eq:eq2}
D_1 ^{(n)} F(\vee _n X) \to F(\vee _n X) \stackrel{F(+)}{\to} F(X).
\end{equation} 
Since the `$+$'-map is $\Sigma_n$-equivariant, the Equation~\ref{eq:eq2} induces a morphism
\begin{equation}
\label{eq:eq3}
D_n F(X) \simeq D_1 ^{(n)} cr_n F(X)_{h\Sigma_n} \stackrel{\simeq}{\to} 
D_1 ^{(n)} F(\vee _n X)_{h\Sigma_n} \to F(X) \to P_{n+1} F(X),
\end{equation} 
where the first equivalence $D_n F(X) \simeq D_1 ^{(n)} cr_n F(X)_{h\Sigma_n}$ is by 
Proposition~\ref{prop:deriv}. We induct on $n$ to show that these maps produce 
the desired splittings. 

For $n=2$, recall that the morphism $D_1 F(X) \to F(X) \to P_2 F(X)$ is a splitting to 
$P_2 F(X) \stackrel{q_2}{\to} D_1 F(X)$, hence $P_2 F(X)$ is equivalent to $D_1 F(X) \vee D_2 F(X)$. 

Now suppose for all $n \leq k$ we have that $P_n F(X)$ is equivalent to the coproduct of 
layers $D_1 F(X) \vee \cdots \vee D_n F(X)$. Consider the diagram 
$$
\xymatrix{
D_{k+1} F(X) \simeq D_1 ^{(k+1)} cr_{k+1} F(X)_{h\Sigma_{k+1}} \ar[r]^-{d_{k+1}}
&
P_{k+1} F(X) \ar[r]^{q_{k+1}} \ar[d]^{q_{k+1}}
&
P_k F(X) \ar[d]^{q_{k}}\\
D_{k} F(X) \simeq D_1 ^{(k)} cr_{k} F(X)_{h\Sigma_{k}} \ar[r]^-{d_k} \ar[ur]
&
P_{k} F(X) \ar[r]^{q_{k}}
&
P_{k-1} F(X). \ar@/^1pc/[l]
}
$$
To see that $P_{k+1} F(X)$ splits we simply need to show that the composite 
\begin{equation}
\label{eq:eq4}
P_k F(X) \simeq P_{k-1} F(X) \vee D_1 ^{(k)} cr_{k} F(X)_{h\Sigma_{k}} \to 
F(X) \stackrel{p_{k+1}}{\to} P_{k+1} F(X) \stackrel{q_{k+1}}{\to} P_k F(X)
\end{equation}
is equivalent to the identity map. Here the first map on the component $P_{k-1} F(X)$ 
exists because by inductive hypothesis, $P_{k-1} F(X)$ is equivalent to the coproduct 
of layers $D_1 \vee \cdots \vee D_{k-1}$, and the maps on layers are defined via 
Equation~\ref{eq:eq3}. It is enough to prove that  Map~\ref{eq:eq4} is equivalent 
to the identity map on the component $D_1 ^{(k)} cr_{k} F(X)_{h\Sigma_{k}}$, because 
the identity on the component  $P_{k-1} F(X)$ follows by inductive hypothesis.
In other words, we need show that the composite
\begin{equation}
\label{eq:eq5}
D_1 ^{(k)} cr_{k} F(X)_{h\Sigma_{k}} \to F(X) \stackrel{p_{k}}{\to} P_k F(X), 
\end{equation}
which is the restriction of Map~\ref{eq:eq4} to $D_1 ^{(k)} cr_{k} F(X)_{h\Sigma_{k}}$,
is equivalent to the map  
$D_1 ^{(k)} cr_{k} F(X)_{h\Sigma_{k}} \stackrel{d_k}{\to} P_k F(X)$ of the Goodwillie
tower of $F(X)$. 
Consider the commutative diagram
\begin{equation}
\label{diag:pcross}
\xymatrix{
cr_k F(X) \ar[r]^f \ar[d]^{p_k}
&
F(X) \ar[d]^{p_k} \\
P_k cr_k  F(X) \ar[r]^{P_k(f)} \ar@/^1pc/[u]^s
&
P_k F(X) 
}
\end{equation}
where $f$ is the `+'-map $cr_k F(X) \to F(X)$,  
two vertical maps pointing down are from Goodwillie towers of $F(X)$ and 
$cr_k F(X)$, and the upward vertical map $s$ exists by Statement {\bf 1} of the 
Theorem since by 
Proposition~\ref{prop:deriv}, $P_k cr_k  F(X)$ is equivalent to $D_1 ^{(k)} cr_{k} F(X)$.
Consequently, $s$ is simply a successive application of sections to the derivative map 
as in Equation~\ref{eq:eq2}, and thus, is itself a section to $p_k$.    

Observe that the composite $p_k \circ f \circ s:D_1 ^{(k)} cr_{k} F(X)\to P_k F(X)$ 
factors through the homotopy orbits to produce the Map~\ref{eq:eq5}, while the map 
$P_k (f)$ induces  $d_k: D_1 ^{(k)} cr_{k} F(X)_{h\Sigma_{k}} \to P_k F(X)$.
Hence, the commutativity of Diagram~\ref{diag:pcross} implies the  desired equivalence 
of the Map~\ref{eq:eq5} and the map $d_k$. 

Thus, we are allowed to conclude that the map $P_k \to P_{k+1}$ given in Equation~\ref{eq:eq4}
is a splitting to $q_{k+1}$, proving that the Goodwillie tower of $F$ at $X$ splits.

Now let $Y=\vee _t X$ for some $t>1$. By what we just proved, to show that the Goodwillie 
tower of $F$ splits at $Y$, it is enough to produce splittings for derivative maps 
$F(\vee _n Y) \to D_1 ^i F(\vee _n Y)$. We will only consider the case $i=1$ and refer to 
the symmetry of arguments for all other $i$'s.

Recall that  $D_1 ^1 F(\vee _n Y)=D_1 ^1 F(Y \vee \vee _{n-1} Y)$ is linear in the first 
variable. Hence,
\begin{equation*}
D_1 ^1 F(Y \vee \vee _{n-1} Y) \simeq D_1 ^1 F(X_1 \vee \vee _{n-1} Y) \vee \cdots
\vee D_1 ^1 F(X_t \vee \vee _{n-1} Y),
\end{equation*}
where $X_1=\cdots =X_t=X$. The indexing is introduced to help keep track of different 
summands. By Statement {\bf 1}, we have a morphism
\begin{equation*}
D_1 ^1 F(X \vee \vee _{n-1} Y)=D_1 ^1 F(X \vee \vee _{t(n-1)}X)\to F(X \vee \vee _{t(n-1)}X),
\end{equation*}
which is a splitting to the derivative map 
$F(X \vee \vee _{t(n-1)}X) \to D_1 ^1 F(X \vee \vee _{t(n-1)}X)$. Thus, we get a morphism
\begin{eqnarray*}
\label{eq:array1}
{D_1 ^1 F(Y \vee \vee _{n-1} Y) \simeq \bigvee_{i=1} ^t D_1 ^1 F(X_i \vee \vee _{n-1} Y) \to
 \bigvee_{i=1} ^t F(X_i \vee \vee _{n-1} Y) = 
 \bigvee_{i=1} ^t F(X_i \vee \vee _{t(n-1)} X)} \\ 
 \to
 F(\bigvee_{i=1} ^t X_i \vee \vee _{t(n-1)} X) = 
F(Y \vee \vee _{n-1} Y), \nonumber
\end{eqnarray*}
where the last map is induced by obvious inclusions 
$X_i \vee \vee _{t(n-1)} X \to \vee_{i=1} ^t X_i \vee \vee _{t(n-1)} X$. 
To see that this map is a section to the derivative 
$F(Y \vee \vee _{n-1} Y) \to D_1 ^1 F(Y \vee \vee _{n-1} Y)$, we need to show that for each 
$i$, the composite 
\begin{eqnarray}
\label{eq:array2}
{D_1 ^1 F(X_i \vee \vee _{n-1} Y) \to F(X_i \vee \vee _{n-1} Y) \to 
F(\bigvee_{i=1} ^t X_i \vee \vee _{n-1} Y) = F(Y \vee \vee _{n-1} Y)}\\
\to D_1 ^1 F(Y \vee \vee _{n-1} Y) \to D_1 ^1 F(X_i \vee \vee _{n-1} Y) \nonumber
\end{eqnarray}
is equivalent to the identity. Observe that the last two maps in the above composite 
fit into the following commutative diagram
$$
\xymatrix{
F(Y \vee \vee _{n-1} Y) \ar[r] \ar[d]
&
D_1 ^1 F(Y \vee \vee _{n-1} Y) \ar[d]\\
F(X_i \vee \vee _{n-1} Y) \ar[r]
&
D_1 ^1 F(X_i \vee \vee _{n-1} Y)
}
$$
where the vertical maps are induced by projections $Y= \bigvee_{i=1} ^t X_i \to X_i $.
Consequently, we can rewrite the Composite~\ref{eq:array2} as 
\begin{eqnarray*}
\label{eq:array3}
{D_1 ^1 F(X_i \vee \vee _{n-1} Y) \to F(X_i \vee \vee _{n-1} Y) \to 
F(\bigvee_{i=1} ^t X_i \vee \vee _{n-1} Y)}\\
\to F(X_i \vee \vee _{n-1} Y) \to D_1 ^1 F(X_i \vee \vee _{n-1} Y),
\end{eqnarray*}
which is the identity since the first map 
$D_1 ^1 F(X_i \vee \vee _{n-1} Y) \to F(X_i \vee \vee _{n-1} Y)$ is a section to the derivative 
map. 
\end{proof}

\section{Algebras over Operads and Bar Construction}
\label{sect:algbar}
We begin this section by recalling the definition of operads. Operads can be defined 
in any symmetric monoidal category. However, since our interests lie primarily in 
the category of $S$-modules (or more generally in the category of $A$-modules), we 
will focus our discussion around these categories, though the category of chain complexes
$Ch(K)$ over a commutative ring $K$ will often get utilized as well to produce examples. 
Recall that the symmetric monoidal operation
in the category of $S$-modules is the smash product $\wedge$.

\begin{definition}
\label{def:operad}
[See~\cite{Kriz}]
An operad is a sequence of objects ${\bf a}(k)$, $k \geq 0$, carrying an action
of symmetric groups $\Sigma_k$, with products
$$\gamma:{\bf a}(k) \wedge {\bf a}(j_1)\wedge \cdots \wedge {\bf a}(j_k) \to
  {\bf a}(j_1 + \cdots + j_k)$$
which are unital, equivariant and associative in the following sense.

(a) The following associativity diagrams commute, where $\Sigma j_s =j$ and 
$\Sigma i_t =i$; also set $g_s = j_1 + \cdots + j_s$ and 
$h_s = i_{g_{s-1} +1}+ \cdots + i_{g_s}$ for $1 \leq s \leq k$:
$$
\xymatrix{
{\a}(k) \wedge (\bigwedge _{s=1} ^k {\a}(j_s)) \wedge (\bigwedge _{r=1} ^j {\a}(i_r))
\ar[rr]^-{\gamma \wedge Id} \ar[dd]_{shuffle}
&&
{\a}(j) \wedge (\bigwedge _{r=1} ^j {\a}(i_r)) \ar[d]^{\gamma}\\
&&
{\a}(i) \\
{\a}(k) \wedge (\bigwedge _{s=1} ^k {\a}(j_s) \wedge (\bigwedge _{q=1} ^{j_s} {\a}(i_{g_{s-1}+q})))
\ar[rr]_-{Id \wedge (\wedge_s \gamma)}
&&
{\a}(k) \wedge (\bigwedge _{s=1} ^k {\a}(h_s)) \ar[u]_{\gamma}
}
$$

(b) The following unit diagrams commute:

\hspace{40pt}
$
\xymatrix{
{\a}(k) \wedge S^{\wedge k} \ar[r]^{\simeq} \ar[d]_{Id \wedge \eta ^k}
&
{\a}(k)\\
{\a}(k) \wedge {\a}(1)^{\wedge k} \ar[ur]_{\gamma} 
&
}
$
\hspace{80pt}
$
\xymatrix{
S \wedge {\a}(j) \ar[r]^{\simeq} \ar[d]_{\eta \wedge Id}
&
{\a}(j)\\
{\a}(1)\wedge {\a}(j) \ar[ur]_{\gamma} 
&
}
$

(c) The following equivariance diagrams commute, where $\sigma \in \Sigma_k$,
$\tau_s \in \Sigma_{j_s}$, the permutation $\sigma(j_1,  \cdots, j_k) \in \Sigma_j$
permutes $k$ blocks of letters as $\sigma$ permutes $k$ letters, and 
$\tau_1 \oplus \cdots \oplus \tau_k \in \Sigma_k$ is the block sum:

$
\xymatrix{
{\a}(k) \wedge (\bigwedge _{s=1} ^k {\a}(j_s)) 
\ar[r]^-{\sigma \wedge \sigma^{-1}} \ar[d]_{\gamma}
&
{\a}(k) \wedge (\bigwedge _{s=1} ^k {\a}(j_{\sigma _s})) \ar[d]^{\gamma}\\
{\a}(j) \ar[r]^{\sigma(j_{\sigma(1)}, \cdots , j_{\sigma (k)})}
&
{\a}(j)
}
$
\hspace{20pt}
$
\xymatrix{
{\a}(k) \wedge (\bigwedge _{s=1} ^k {\a}(j_s)) 
\ar[r]^-{Id \wedge \wedge _1 ^k \tau _i} \ar[d]_{\gamma}
&  
{\a}(k) \wedge (\bigwedge _{s=1} ^k {\a}(j_s)) \ar[d]^{\gamma}\\
{\a}(j) \ar[r]^{\tau_1 \oplus \cdots \oplus \tau_k }
&
{\a}(j)
}
$

\end{definition}

Our main interest is in the categories of algebras over operads. For detailed good 
discussions on general theory of algebras over operads we refer to~\cite{Kriz} 
and~\cite{Getz}. Here we recall one of the equivalent definitions with a slight 
modification. In what follows, we are mostly going to be concerned with non-unital 
algebras. This corresponds to the additional hypothesis of ${\a}(0) \cong 0$ in the above 
definition, which will be assumed throughout this paper unless specified otherwise.   

We denote the triple associated to $\a$ by $T_{\a}$:
$$T_{\a} (X) = \bigoplus _{n=0} ^{\infty} {\a} (n) \wedge _{\Sigma_n} X^{\wedge n}.$$
Of course, as pointed out, for us the $0-th$ summand is redundant. 
Let $\C _{\a}$ be the category of algebras over the operad $\a$, which coincides with 
the category of algebras over the triple $T_{\a}$. In other words, if $B$ is in $\C_{\a}$, 
it is equipped with a collection of maps
$$ {\a} (n) \wedge _{\Sigma_n} B^{\wedge n} \to B,$$
which make the obvious coherence diagrams commute. We will refer to algebras over an 
operad $\a$ as $\a$-algebras. 

To get a better homotopical control over the category of algebras $\C _{\a}$, we 
consider only the operads ${\a}$ which arise from operads in simplicial sets. Before
making this statement more explicit we remark that this is not very restrictive, as
nearly all operads in the category of spectra that appear in literature are induced 
by an operad in simplicial sets. 

Observe that if ${\bf b}= \{{\bf b}(n)\}$ is an operad in the category of simplicial 
sets, then ${\a}$ with ${\a}(n) = S \otimes {\bf b}(n)$ is an operad in $\M_S$. 
Here we employed the fact that the category of $S$-modules is a tensored category. 
To see the operadic multiplications, recall that there are natural isomorphisms 
$$(S \otimes {\bf b}(k)) \wedge (S \otimes {\bf b}(j_1)) \wedge \cdots \wedge 
(S \otimes {\bf b}(j_k)) \cong S \otimes {\bf b}(k) \otimes {\bf b}(j_1) 
\otimes \cdots \otimes {\bf b}(j_k).$$
Hence we can use the multiplication of ${\bf b}$ to produce the operadic multiplication 
on ${\a}$. This is what we mean when discussing operads in $\M_S$ induced by  
operads in the category of simplicial sets.    

From this point on, unless specifically said otherwise, the word ``operad'' means 
a cofibrant operad in the category of $S$-modules that arises from simplicial sets.
 
The following Proposition is part of Proposition 1.6 of~\cite{Goerss}.

\begin{proposition}
\label{prop:modcat}
If we define a map $X \to Y$ of ${\a}$-algebras to be a weak equivalence or 
a fibration if it is a weak equivalence or a fibration of $S$-modules, 
then the category of ${\a}$-algebras $\C_{\a}$ becomes a cofibrantly 
generated simplicial model category.    
In particular, it is equipped with a cofibrant replacement functor.
\end{proposition}

The proof of this proposition is sketched in Section 1 of~\cite{Goerss}, 
for more details see also~\cite{Berger}.

The following is a variant of the two sided bar construction of P.May (~\cite{May}).

Denote by $\mu_{\a}$ and $\eta_{\a}$ the multiplication and the unit respectively of 
the cotriple $T_{\a}$. When no confusion can arise we will omit the subscript 
$\a$ from the notation.  
For an $\a$-algebra $(C, \xi)$, let $B_{\ast} ^{\a} C$ denote the simplicial object with 
$n$'th term $B_{n} ^{\a} C = T_{\a}^{n+1} C$ and face and degeneracy operators given 
by 

$d_i = T_{\a}^{i} \mu_ {T_{\a}^{n-i-1}}$, for  $0 \leq i < n$, and $d_n = T_{\a}^{n} \xi$

$s_i = T_{\a}^{i+1} \eta_ {T_{\a}^{n-i}}$, for $0 \leq i \leq n$.

This construction can be interpreted in the language of~\cite{May} if we consider 
$T_{\a}: \M_{A} \to \C_{\a}$ as the functor left adjoint to the forgetful functor 
$U:\C_{\a} \to \M_{A}$. Then $B_{\ast} ^{\a} C$ is the same as the object 
$B_{\ast} ( T_{\a}, U T_{\a}, UC)$ in notation of~\cite{May}.
Consequently, we can consider $\xi: B_{\ast} ^{\a} C \to C$ as an augmented simplicial 
$A$-module for which one defines a contraction using the unit map 
$\eta: C \to T_{\a} C$. Here we have omitted the forgetful functor $U$ from our notations 
and will continue to do so, as it is evident from the context. 

Of course for the bar construction 
$B_{\ast} ^{\a} C$ to be computationally useful we need to be allowed to work with it 
level-wise. To make this more precise we recall the notion of properness. To avoid confusion,
here (and in the future) we use the term `c-cofibration' (for classical cofibration)
to distinguish it from cofibrations that are part of the model category structure. 
In other words, a c-cofibration of $A$-modules is simply a map $i: M \to \overline{M}$ of 
$A$-modules that satisfies the homotopy extension property in the category of 
$A$-modules. Of course, all cofibrations are c-cofibrations, but not conversely. 

\begin{definition}
Let $K_\ast$ be a simplicial spectrum and let $sK_q \subset K_q$ be the `union' of the 
subspectra $s_jK_{q-1}$, $0 \leq j< q$. A simplicial $A$-module $K_\ast$ is proper 
if the canonical map of $A$-modules $sK_q \to K_q$ is a c-cofibration for each 
$q \geq 0$.

\end{definition}

The main reason that proper simplicial $A$-modules $K_\ast$ are computationally 
useful is that one is allowed to work with them level-wise. More precisely,
one can use the simplicial filtration to construct a well-behaved 
spectral sequence that converges to $\pi_\ast (E \wedge K)$ for any spectrum $E$
(Theorem X.2.9 of ~\cite{EKMM}). In particular, if $f:K_\ast \to L_\ast$ is a map
of proper simplicial $A$-modules which is a weak equivalence level-wise, then 
the geometric realization $|f|$ of $f$ is also a weak equivalence. 

\begin{lemma}
For a cofibrant $\a$-algebra $C$, $B_\ast ^{\a} C$ is a proper simplicial $A$-module.
\end{lemma}

\begin{proof}
The condition of properness involves only the degeneracy operators (and not
the face maps) of a simplicial $A$-module. As it is evident from 
the definition of the bar construction the degeneracies are constructed 
from the unit map $\eta : C \to T_{\a} C$. This map has an obvious section 
$T_{\a} C \to C$ (given by the multiplication map) which is a map of $A$-modules. 
Thus, it satisfies the homotopy 
extension property. Similarly one shows that all degeneracies are 
c-cofibrations.  
\end{proof}

Note that the notion of geometric realization of the 
simplicial $\a$-algebra $B_\ast ^{\a} C$ is somewhat ambiguous, as one could realize this object
in the category of $A$-modules, or alternatively, in the category of $\a$-algebras, i.e. 
internally. We denote the realization in the category of $\a$-algebras by $|-|_{\a}$. 
The following proposition shows that the two realizations are the same. 
The proof presented here is based on an argument suggested by Paul Goerss. 
Proposition VII.3.3 of~\cite{EKMM} provides a different, categorical approach.  

\begin{proposition}
The geometric realization in the category of $\a$-algebras is isomorphic to the 
geometric realization in the category of $A$-modules.
\end{proposition}
   
\begin{proof}
First observe that it is sufficient to prove that the natural map 
\begin{equation}
\label{eq:geomalg}
|T_{\a} M_\ast | \to T_{\a} |M_\ast |
\end{equation}
is an isomorphism, where $M_\ast$ is a simplicial $A$-module. Indeed, if this is the case then
for any simplicial $\a$-algebra $C_\ast$, the geometric realization as modules is an 
$\a$-algebra via the composite (multiplication) map
$$T_{\a} |C_\ast | \simeq |T_{\a} C_\ast | \to |C_\ast |.$$ 
Consequently, by an argument similar to that in the proof of Proposition 3.3 of~\cite{EKMM}, we
get that for any simplicial $\a$-algebra $D_\ast$, there is an isomorphism
$${\C}_{\a} (|C_\ast |, D_\ast) \cong s{\C} (C_\ast, F_{\a} (\triangle _{\ast}, D_\ast)),$$
where  $F_{\a}$ is the internal function space functor. Recalling that 
 $F_{\a}(\triangle _{\ast},-)$ is right adjoint to the geometric realization functor 
$|-|_{\a}$, we conclude that, by Yoneda's lemma, $|C_{\ast}|_{\a} \simeq |C_{\ast}|$. 

Thus, it remains to show that Morphism~\eqref{eq:geomalg} is an equivalence. Note that 
as a left adjoint, geometric realization commutes with coproducts and the functor 
${\a}(n) \wedge _{\Sigma_n}(-)$. Hence it is enough to show that 
$|M_{\ast} ^{\wedge n}|$ is equivalent to $|M_{\ast}| ^{\wedge n}$. Further recalling that 
geometric realization commutes with the spectrification functor, we reduce the problem 
to showing that for any pointed simplicial sets $X_\ast$ and $Y_\ast$, the natural map 
$|X_\ast \wedge Y_\ast| \to |X_\ast | \wedge |Y_\ast|$ is an equivalence. 

However, this is an immediate consequence of the Eilenberg-Zilber theorem, once we observe that 
both $|X_\ast \wedge Y_\ast|$ and $|X_\ast | \wedge |Y_\ast|$ can be obtained from the 
bisimplicial space
$$([m], [n]) \mapsto X_m \wedge Y_n,$$
with $|X_\ast \wedge Y_\ast|$ being the geometric realization of the diagonal, and 
$|X_\ast | \wedge |Y_\ast|$ resulting from realizing in the $m$-direction first and then
in the $n$-direction.  

\end{proof}

As an immediate consequence of this proposition, we have that $|B_{\ast} ^{\a} C|$ and 
$C$ are weakly equivalent as $\a$-algebras. Indeed, we already observed that  $|B_{\ast} ^{\a} C|$ 
and $C$ are homotopy equivalent as $A$-modules (via the contraction to the augmentation map),
hence the weak equivalence as algebras follows from the above proposition.

\section{The Forgetful Functor $U_{\a}: C_{\a} \to \M_A$}
\label{sect:forget}

In this section we construct a tower of functors that approximates the forgetful
functor $U_{\a}: C_{\a} \to \M_A$ from the category of $\a$-algebras to the category 
of $A$-modules. As noted in the previous section, we are still assuming that our operads arise 
from simplicial sets. 

Roughly speaking, we say that the functor $F(X)$ is approximated 
by the tower of functors $P_n F(X)$ if they assemble into a diagram
$$
\xymatrix{
&
\cdots \ar[d]^{q_{n+1}}\\
F(X) \ar[r]^{p_n} \ar[rd]^{p_{n-1}}
&
P_n F(X) \ar[d]^{q_{n}}\\
&
P_{n-1} F(X) \ar[d]\\
&
\cdots
}$$ 
where $p_i$'s and $q_i$'s are natural transformations, such that under suitable conditions
the homotopy limit of $\{P_n F(X) \}$ is weakly equivalent to $F(X)$. 

As noted, one of the objectives of this work is to determine when the 
given $\a$-algebra decomposes into a direct sum, or more precisely, what are sufficient
conditions for the algebra to be free. The idea is to construct a tower of functors 
approximating the above forgetful functor, and search for conditions under which this 
tower splits. We will return to the relationship between freeness and splitting towers,
as well as to the precise meaning of the term `approximate' at a later point.  

The motivation for the `approximating' tower of functors comes from Goodwillie Calculus. 
In particular, a similar tower is constructed in ~\cite{Randy} for the discrete commutative case
and in ~\cite{Minas} for the $E_\infty$ case. However, as noted, the primary objective 
of this part of the paper is to develop splitting criteria for algebras over operads, and in a 
large number of cases we are able to obtain and give a nice description of such splittings
without having to resort to the terminology of Goodwillie Calculus. 

The downside is an additional condition on the operad. For algebras over a general operad, we are 
still able to prove a splitting result, but the pieces into which our algebras decompose 
are somewhat less transparent, and Calculus language could not be avoided. It should be noted 
though that these `two' decompositions are essentially the same in their common range; we will 
comment more on this at a later stage. 

On an organizational point, this and the section following it, are constructed in such a manner
that the results involving Goodwillie Calculus are clearly marked and can be omitted 
by the reader (if he/she chooses to do so) without affecting the rest of these two sections.

We begin by defining a functor $Q_n:\C_{\a} \to \M_A$ for all $n>1$ by the following 
pushout diagram in the category of $A$-modules

$$
\xymatrix{
{\a}(n) \wedge_{\Sigma_n} C^{\wedge n} \ar[r] \ar[d]^{\mu}
&
\ast \ar[d]\\
C \ar[r] 
&
Q_n(C) .
}
$$
This functor already has some of the ingredients we require, namely it is equipped with 
a natural transformation $U_{\a} \to Q_n$, which is evident from the definition of $Q_n$. 
Furthermore, under suitable connectivity assumptions on an $\a$-algebra $C$, the 
connectivity of the maps $U_{\a}(C) \to Q_n(C)$ increases with $n$. However, we
want to work with $C_{\a}$ not only up to isomorphism, but up to a weak equivalence. Thus, 
we desire functors which preserve weak equivalences. Usually this is achieved 
by considering the derived version of the functor, i.e. functors get evaluated not at 
the objects themselves but at their cofibrant replacements. In addition, we would like our 
functors to be computationally friendly, which is often accomplished by 
constructing a cofibrant replacement via a resolution by free objects.

To that end, we introduce the following definition. In what follows, $\Gamma C$ is the cofibrant 
replacement functor in the category of $\a$-algebras, and it exists by Proposition~\ref{prop:modcat}.
\begin{definition}
For an $\a$-algebra $C$, define $I/I^n(C)$ to be the $A$-module
$$I/I^n(C) = Q_n(\Gamma B^{\a} _{\ast} \Gamma C)$$
where $Q_n$ is applied level-wise. We will also write  

\hspace{100pt} $I^n/I^{n+1}(C) \stackrel{def}{=} hofiber[I/I^{n+1}(C) \rightarrow I/I^n(C)]$.
\end{definition}

\begin{remark}
Note that even when $C$ is cofibrant as an $A$-module, $B^{\a} _{\ast} C$ 
is not necessarily cofibrant. Of course, $T_{\a}$ takes cofibrant $A$-modules to cofibrant 
$\a$-algebras, however, cofibrant $\a$-algebras are not necessarily cofibrant as $A$-modules. 
Consequently, when we apply $T_{\a}$ again, the resulting $\a$-algebra may no longer be 
cofibrant. Hence, since simplices of $B^{\a} _{\ast} C$ are formed by multiple 
applications of $T_{\a}$, they may not be cofibrant. Thus, we are forced to replace 
$B^{\a} _{\ast} \Gamma C$ by  $\Gamma B^{\a} _{\ast} \Gamma C$, where the outer $\Gamma$ is 
applied level-wise. Naturally, in doing so, we lose some of the computational advantage.
However, as is illustrated in the proof of the next proposition, this problem can be overcome.  
\end{remark}

The following Proposition~\ref{prop:fiber} gives us a good description of the layers (fibers) of the 
approximating tower of functors. More precisely it allows us to express the higher layers
in terms of the first one. As a condition on the operad is required, we introduce a definition first. 

\begin{definition}
\label{def:prim}
We say that the operad ${\a}$ is primitively generated if the square
$$
\xymatrix{
{\a}(n) \wedge_{\Sigma_n} T_{\a}(X) ^{\wedge n} \ar[r] \ar[d]
&
T_{\a}(X) \ar[d]\\
\ast \ar[r]
&
\bigoplus_{i \leq n-1} {\a}(i) \wedge_{\Sigma_i} X^{\wedge i}
}
$$ 
is Cartesian for all $X$ and $n \geq 1$. 
\end{definition}  

Note in particular, that for primitively generated operads, the multiplication map 
${\a}(n) \wedge {\a}(1) \wedge \cdots \wedge {\a}(1) \to {\a}(n)$ is an equivalence. 

\begin{remark}
This is both a remark and an example as its purpose is to explain our 
choice of terminology in Definition~\ref{def:prim}, as well as to show that  
quadratic operads (see~\cite{GinKap}) are primitively generated.
To do so, we analyze the 
primitively generated operads in the category of differential graded complexes 
more closely. Let $\a$ be such an operad. To understand the top horizontal map in the 
above diagram better, observe that 
$$(T_{\a} X)^{\wedge n} = \bigoplus_{(j_1, \cdots , j_n)} {\a}(j_1) \wedge \cdots \wedge {\a}(j_n)
\wedge_{\Sigma_{j_1} \times \cdots \times \Sigma_{j_n}} C^{\wedge (j_1 + \cdots + j_n)}.$$
Consequently, for say $n=2$, that map is given by components
$${\a}(2) \wedge [{\a}(s) \wedge {\a}(t) \oplus {\a}(t) \wedge {\a}(s)]
\wedge_{\Sigma_2 \ltimes (\Sigma_{t} \times \Sigma_{s})} X^{\wedge (s+t)} \to 
{\a}(s+t) \wedge_{\Sigma_{s+t}} X^{\wedge (s+t)},$$
which are a combination of the `$+$'-map and the operadic multiplication. Since for all $X$
all summands of degree $n$ and higher must get annihilated in the above diagram, 
and for $n=2$, the component ${\a}(k) \wedge_{\Sigma_k} X^{\wedge k}$ gets ``hit'' by 
$$ \bigoplus_{t+s=k} {\a}(2) \wedge [{\a}(s) \wedge {\a}(t) \oplus {\a}(t) \wedge {\a}(s)]
\wedge_{\Sigma_2 \ltimes (\Sigma_{t} \times \Sigma_{s})} X^{\wedge (s+t)},$$
with $t$ and $s$ greater than $0$, then we can assert that ${\a}(k)$ is produced by the 
objects ${\a}(t)$ of degree lower than $k$. Inducting down on $k$, we conclude that 
all objects in the sequence of the operad $\a$ are ``produced'' by ${\a}(1)$ and ${\a}(2)$. 
To make this more precise we recall the notion of a free operad on a symmetric sequence 
${\bf v}(l)$. In other words, for each $l$, ${\bf v}(l)$ is an object in the symmetric 
monoidal category with an action of the symmetric group $\Sigma_l$. The forgetful functor from 
the category of operads to symmetric sequences has a left adjoint $\mathbb T$ 
endowed with a natural map $\eta: {\bf v} \to {\mathbb T} {\bf v}$, which satisfies 
the following universal property. If ${\bf p}$ is an operad, and $\bar{\eta}:{\bf v} \to {\bf p}$ 
is a map of symmetric sequences, than there is a unique morphism of operads 
$\phi:{\mathbb T} {\bf v} \to {\bf p}$, such that $\bar{\eta}=\phi \eta$. The operad 
${\mathbb T} {\bf v}$ is referred to as the free operad on ${\bf v}$. A detailed description 
of the free operad functors in terms of trees is provided in~\cite{GinKap}. 
 
In addition, it is possible to define the notion of an operad {\it ideal}, in such a way 
that ${\bf c} \subset {\a}$ is an operad ideal if and only if the operadic multiplication 
of $\a$ induces an operadic multiplication on the quotient symmetric sequence ${\a}/{\bf c}$. 

Thus, using this new language, we assert that the primitively generated operads are 
quotients of the free operad on symmetric sequences  ${\bf v}(l)$ with ${\bf v}(l)=0$ 
for all $l>2$. Conversely, suppose the operad $\a$ is such that the operadic multiplication
${\a}(2) \wedge {\a}(1) \wedge {\a}(1) \to {\a}(2)$ is an equivalence and 
the generators of the operad are all in degrees $1$ and $2$.
Then the top horizontal map 
$${\a}(n) \wedge_{\Sigma_n} T_{\a}(X)^{\wedge n} \to   T_{\a}(X)$$
of the diagram in Definition~\ref{def:prim} misses the components of degree less than $n$,
and ``hits'' all the summands starting from $n$, thus making the diagram Cartesian. 
 
{\it Quadratic} operads provide examples of these as they are defined to 
be the operads whose generators are concentrated in degree $2$ while the relations are in 
degree $3$; see~\cite{GinKap}. Consequently, all the well known quadratic operads, such the 
{\it associative}, {\it commutative}, and {\it Lie} operads are examples of primitively 
generated operads.

\end{remark}

Now we state the promised result. 

\begin{proposition}
\label{prop:fiber}
Let ${\a}$ be a primitively generated operad. 
Then for every $\a$-algebra $C$, there is a weak equivalence of simplicial objects

\hspace{140pt}
$I^n/I^{n+1}(C) \simeq {\a}(n) \wedge_{h\Sigma_n} \bigwedge ^n I/I^2(C),$

\noindent
where $\bigwedge ^n I/I^2$ is the $n$'th smash power of $I/I^2$.
\end{proposition}

This is a generalization of Proposition 2.4 of~\cite{Minas}, where a similar result is proved 
for $E_\infty$ operads. In proving this proposition we use some of the techniques developed in 
~\cite{Minas} and~\cite{Maria}. We will need a key technical lemma, which we will prove later, in 
Section~\ref{sect:last}, in order not to distract from our  task at hand. 

\begin{lemma}
\label{lem:tech}
Let $C$ be a cofibrant $\a$-algebra and $\gamma:Y \to T_{\a}^n C$ a cell $A$-module approximation.

{\bf 1.}
Then the induced map 
$${\a}(i) \wedge_{\Sigma_i} Y^{\wedge i} \to {\a}(i) \wedge_{\Sigma_i} (T_{\a} ^n C)^{\wedge i} $$ 
is a weak equivalence for all $i>0$.

{\bf 2.}
The projection 
${\a}(i)  \wedge_{h\Sigma_i} (T_{\a} ^n C)^{\wedge i} \to 
 {\a}(i) \wedge_{\Sigma_i} (T_{\a} ^n C)^{\wedge i}$
is an equivalence for all $n$ and $i$. 
\end{lemma}

\begin{proof} 
(of Proposition~\ref{prop:fiber})
As it was noted earlier, by applying the functor $B_\ast ^{\a}$ we gained some computational
advantages (as it is the analogue of taking a free resolution in discrete algebra), however
some of it was lost when we were forced to take a cofibrant replacement. Hence, the first 
objective is to show that we get equivalent constructions even if we forgo taking cofibrant 
replacements, in other words, we begin by showing that $I/I^n(C)$ is weakly equivalent to 
$Q_n(B_\ast ^{\a} \Gamma C)$ as modules. 

Let $\gamma: Y \to T_{\a} ^k C$ be a cell $A$-module approximation of $T_{\a} ^k C$, where
$C$ is a cofibrant $\a$-algebra. Then by Lemma~\ref{lem:tech}, we get that 
$T_{\a} \gamma: T_{\a} Y \rightarrow T_{\a} T_{\a} ^k C$ is a weak equivalence, since it 
is a coproduct of weak equivalences.  
Hence so is
$\Gamma T_{\a} Y \rightarrow \Gamma T_{\a} T_{\a} ^k C$. In fact, since the two 
algebras involved are cofibrant, this last map is a simplicial homotopy equivalence.  
Consequently, the map $Q_n(\Gamma T_{\a} Y) \rightarrow Q_n(\Gamma T_{\a} T_{\a} ^k C)$
is also a simplicial homotopy equivalence, as any functor preserves simplicial homotopy equivalences. 

Now consider the following commutative diagram:

$$
\xymatrix{
Q_n(\Gamma T_{\a} Y) \ar[r] \ar[d]
&
Q_n(\Gamma T_{\a} T_{\a} ^k C) \ar[d] \\
Q_n(T_{\a} Y) \ar[r]
&
Q_n( T_{\a} T_{\a} ^k C).
}
$$
As we just argued, the top horizontal map is a homotopy equivalence. Since $T_{\a} Y$ is a 
cofibrant $\a$-algebra, the cofibrant replacement weak equivalence  
$\Gamma T_{\a} Y \to T_{\a} Y$ is in fact a homotopy equivalence. Hence the left vertical 
map is also a homotopy equivalence. To see that the bottom horizontal map is a weak 
equivalence, we analyze $Q_n(T_{\a} Y)$ more closely. By definition, it is the cofiber 
of the map 
$${\a}(n) \wedge_{\Sigma_n}
 [\bigoplus _{i=1}^{\infty} {\a}(i) \wedge_{\Sigma_i} Y^{\wedge i}]^{\wedge n} \to
 \bigoplus _{i=1}^{\infty} {\a}(i) \wedge_{\Sigma_i} Y^{\wedge i},$$  
which, since $\a$ is primitively generated, is equivalent to 
$\bigoplus_{i \leq n-1} {\a}(i) \wedge_{\Sigma_i} Y^{\wedge i}$.
Similarly,
$$Q_n( T_{\a} T_{\a} ^k C) \simeq 
\bigoplus_{i \leq n-1} {\a}(i) \wedge_{\Sigma_i} (T_{\a} ^k C) ^{\wedge i}.$$
The equivalence of the bottom horizontal arrow follows from Lemma~\ref{lem:tech}. 
Given that three of the arrows in the above diagram are equivalences, we conclude that 
the fourth one, $Q_n( \Gamma T_{\a} T_{\a} ^k C) \to Q_n( T_{\a} T_{\a} ^k C)$, is also a 
weak equivalence, which proves that $I/I^n (C)$ and $Q_n(B_\ast ^{\a} \Gamma C)$ have 
equivalent simplices. We conclude that the map $I/I^{n+1} (C) \to I/I^n (C)$ is a weakly 
equivalent to a fibration which, on the level of $(k+1)$-simplices is given by projections
$$
{\a(1)} \wedge T_{\a} ^k (C) \vee \cdots \vee 
{\a(n)} \wedge_{\Sigma_{n}} (T_{\a} ^k (C))^{\wedge n} \to
{\a(1)} \wedge T_{\a} ^k (C) \vee \cdots \vee 
{\a(n-1)} \wedge_{\Sigma_{n-1}} (T_{\a} ^k (C))^{\wedge n-1}.$$
Thus, for all $n \geq 1$, $I^n/I^{n+1} (C)$ is equivalent to a simplicial $A$-module with 
${\a}(n) \wedge_{\Sigma_{n}} (T_{\a} ^k (C))^{\wedge n}$ for $(k+1)$-simplices. In particular,
the simplices of $I/I^2$ are ${\a(1)} \wedge T_{\a} ^k (C)$. Furthermore, note that by 
Lemma~\ref{lem:tech}, the orbits and homotopy orbits with respect to symmetric groups 
$\Sigma_n$ of the objects involved (i.e. of operadic powers of  $T_{\a} ^k (C)$) are equivalent.
Consequently, given that $\a$ is primitively generated and hence the multiplication maps 
${\a}(n) \wedge {\a}(1) \wedge \cdots \wedge {\a}(1) \to {\a}(n)$ are equivalences
(see the comments following Definition~\ref{def:prim}), we get that 
$I^n/I^{n+1} (C) \simeq {\a}(n) \wedge_{h\Sigma_n} \bigwedge^n (I/I^2(C))$
as desired.

\end{proof}

We return to the question of `approximating' towers. The following 
proposition makes our previously used language precise. Roughly speaking, it states that 
under favorable conditions, the tower of functors $\{I/I^n\}$ converges to the forgetful 
functor $U_{\a}$. 

\begin{proposition}
\label{prop:conn}
Let $A$ be a connective commutative $S$-algebra, and $\a$ a primitively generated
operad in the category of 
$A$-modules. Then for every $0-connected$ $\a$-algebra $C$, the natural map
$$\phi: U_{\a} (\Gamma C) \simeq U_{\a} (\Gamma B_{\ast} ^{\a} \Gamma C) \to \holim_n I/I^n (C)$$
is a weak equivalence. 
\end{proposition}

\begin{proof}
We begin by observing that $(\Gamma B_{\ast} ^{\a} \Gamma C)^{\wedge n}$ is at 
least $n-1-connected$. 
Indeed, since $A$ is connective, the Cellular Approximation Theorem [Chapter 3 of ~\cite{EKMM}]
allows us to functorially replace $\Gamma B_{\ast} ^{\a} \Gamma C$ by a weakly equivalent
CW $A$-module $M$ with no cells in dimensions below one. 
By properties listed in Section~\ref{sect:last}, the derived 
smash powers of $\Gamma B_{\ast} ^{\a} \Gamma C$ are defined on the point set level. Hence, 
we have that $(\Gamma B_{\ast} ^{\a} \Gamma C)^{\wedge n}$
is equivalent to $M^{\wedge n}$. (Once again, this last statement is not immediate since 
$\Gamma B_{\ast} ^{\a} \Gamma C$ is not necessarily cofibrant as an $A$-module.) Observe that
$M^{\wedge n}$ (and consequently $(\Gamma B_{\ast} ^{\a} \Gamma C)^{\wedge n}$) is 
$(n-1)$-connected, since $M$ has no cells in dimensions below $1$, and thus 
$M^{\wedge n}$ has no cells in dimensions less than $n$. 

Hence, we are allowed to conclude that the connectivity of the maps 
$U_{\a}(\Gamma B_{\ast} ^{\a} \Gamma C) \to I/I^n (C)$ increases with $n$, producing a 
Mittag-Leffler system (see~\cite{Bous}), and thus implying the claim of the proposition. 
\end{proof}

As we will see in the next section, Proposition~\ref{prop:fiber} is a critical ingredient
for obtaining a splitting result for $\a$-algebras. Of course, in that proposition, the  
assumption that the operad $\a$ be primitively generated is also used in an essential way. In fact,
without that assumption, the tower $I/I^n$ that we constructed to approximate the forgetful 
functor $U_{\a}$ from the category of $\a$-algebras to $A$-modules does not enjoy the 
properties required to develop our theory. 

However, as we will show momentarily, the Goodwillie tower of the forgetful functor does. 
In fact, for a primitively generated operad $\a$, the tower $\{ I/I^n \}$ {\it is} the 
Goodwillie tower of the forgetful functor. In other words, the condition on the operad 
simply allowed us to provide a nice description of the Goodwillie tower in terms of 
functors $I/I^n$. 

The general version of our result still requires a condition on the operad, however it is 
very mild as nearly all naturally occurring operads satisfy it. 

\begin{proposition}
\label{prop:genfiber}
Let the operad $\a$ be such that for all $n>0$ the operadic multiplication maps 
$${\a}(n) \wedge {\a}(1) \wedge \cdots \wedge {\a}(1) \to {\a}(n)$$
are equivalences. Then for every $\a$-algebra $C$, the $n$'th layer $D_n$ of the 
Goodwillie tower of the forgetful functor $U_{\a}$ can be described as 
$$D_n U_{\a}(C) \simeq {\a} (n) \wedge_{h\Sigma_n} (D_1 U_{\a}(C))^{\wedge n}.$$ 
\end{proposition}

\begin{proof}
We employ the same strategy as when proving Proposition~\ref{prop:fiber}. In other words,
instead of the Goodwillie tower $P_n U_{\a} (\Gamma C)$ consider the equivalent tower
$P_n U_{\a} (\Gamma B^{\a} _{\ast} \Gamma C)$. As before, we can strip the outer cofibrant 
replacement functor $\Gamma$. Indeed, let $\gamma :Y \to T_{\a} ^k \Gamma C$ be a cell 
$A$-module approximation and consider the commutative square

$$
\xymatrix{
P_n U_{\a} (\Gamma T_{\a} Y) \ar[r] \ar[d]
&
P_n U_{\a} (\Gamma T_{\a} T_{\a} ^k \Gamma C) \ar[d]\\
P_n U_{\a} (T_{\a} Y) \ar[r]
&
P_n U_{\a} (T_{\a} T_{\a} ^k \Gamma C).
}
$$
The left vertical and the top horizontal arrows are equivalences by the same argument as in 
Proposition~\ref{prop:fiber}. To see that the bottom horizontal map is an equivalence, note that
since $T_{\a}$ is coproduct preserving, $P_n U_{\a} (T_{\a} Y) \simeq P_n (U_{\a} T_{\a}) (Y)$. 
Consequently,

$$
P_n U_{\a} (T_{\a} Y) = 
P_n [U_{\a} (\bigoplus_{i=1} ^{\infty} {\a}(i) \wedge_{\Sigma_i} Y^{\wedge i})] \simeq 
P_n [U_{\a} (\bigoplus_{i=1} ^{\infty} {\a}(i) \wedge_{h\Sigma_i} Y^{\wedge i})] \simeq
\bigoplus_{i=1} ^{\infty} P_n ({\a}(i) \wedge_{h\Sigma_i} Y^{\wedge i}).
$$
Here we have used the fact that $Y$ is a cell $A$-module (and thus orbits are equivalent 
to homotopy orbits in the above setup), and that the $n$'th Taylor polynomial 
$P_n$ commutes with coproducts of functors. Now note that 
${\a}(i) \wedge_{h\Sigma_i} Y^{\wedge i}$ is an 
$n$-homogeneous functor, i.e. $P_n({\a}(i) \wedge_{h\Sigma_i} Y^{\wedge i}) \simeq \ast$ if $n<i$,
and $P_n({\a}(i) \wedge_{h\Sigma_i} Y^{\wedge i}) \simeq {\a}(i) \wedge_{h\Sigma_i} Y^{\wedge i}$
for $n \geq i$. Thus,
$$
P_n U_{\a} (T_{\a} Y) \simeq 
\bigoplus_{i=1} ^{\infty} P_n U_{\a}({\a}(i) \wedge_{h\Sigma_i} Y^{\wedge i}) \simeq 
\bigoplus_{i=1} ^n {\a}(i) \wedge_{h\Sigma_i} Y^{\wedge i}. 
$$
Similarly,
$$
P_n U_{\a} (T_{\a} T_{\a} ^k \Gamma C) \simeq 
\bigoplus_{i=1} ^n {\a}(i) \wedge_{h\Sigma_i} (T_{\a} ^k \Gamma C)^{\wedge i},
$$
except here since $T_{\a} ^k \Gamma C$ is not necessarily a cell $A$-module, we make use of 
Lemma~\ref{lem:tech} and properties listed in Section~\ref{sect:last}. 
Now Lemma~\ref{lem:tech} implies that the bottom horizontal arrow of the diagram 
is a weak equivalence. Hence, so is the right vertical map, thus proving that we may strip away
the outer functor functor $\Gamma$. Consequently, it suffices to prove the proposition for 
free $\a$-algebras only. The computations for the free case though are identical to those performed
at the end of the proof of Proposition~\ref{prop:fiber}. 

\end{proof}

We conclude this section by pointing out that Proposition~\ref{prop:conn} has its obvious 
analogue in this generalized set up as well. 

\begin{proposition}
\label{prop:genconn}
Let $A$ be a connective commutative $S$-algebra, and $\a$ an operad in the category of 
$A$-modules, satisfying the condition of Proposition~\ref{prop:genfiber}.
Then for every $0-connected$ $\a$-algebra $C$, the natural map
$$\phi: U_{\a} (\Gamma C) \to \holim_n P_n U_{\a}(\Gamma C)$$
is a weak equivalence. 
\end{proposition}

\begin{proof}
The proof of Proposition~\ref{prop:conn} works here verbatim if we replace 
$I/I^{n+1}(C)$ by $P_n U_{\a}(\Gamma C)$.
\end{proof}

\section{Algebras over Operads and Splitting}
\label{sect:splitalg}

\subsection{On Splitting the Forgetful Functor}

In this subsection we address the question of detecting free objects in the category of 
$\a$-algebras. More precisely, given an $\a$-algebra $C$, when is it equivalent to 
$T_{\a} X$ for some $A$-module $X$? Here, as before, $\a$ is an operad arising from 
simplicial sets.  

\begin{theorem}
\label{theor:algsplit}
Let the operad $\a$ be primitively generated. As before, let 
$U_{\a}: \C_{\a} \to \M_A$ be the forgetful functor from the category of $\a$-algebras 
to the category of $A$-modules. Furthermore, suppose the $\a$-algebra $C$ is such that 
the natural map $U_{\a}(C) \to  I/I^2 (C)$ has a section 
$\phi:I/I^2 (C) \to U_{\a}(C) $ in the category of $A$-modules. Then we have a weak
equivalence of $\a$-algebras
$$\holim_n I/I^n (C) \simeq \prod {\a}(n) \wedge_{\Sigma_n} [I/I^2 (C)]^{\wedge n}.$$
\end{theorem}

\begin{proof}
Note that by applying the functor $T_{\a}$ to the section $\phi:I/I^2 (C) \to U_{\a}(C)$
and composing the resulting $\a$-algebra map with the multiplication map $\mu$ with which the 
algebra $C$ is equipped, we get a composite map $\alpha$ of $\a$-algebras
$$\alpha: T_{\a} [I/I^2 (C)] \stackrel{ T_{\a} \phi}{\longrightarrow} T_{\a}[U_{\a}(C)] 
\stackrel{\mu}{\longrightarrow} C.$$
Consider the map $I/I^n (\alpha)$:
$$I/I^n(T_{\a} [I/I^2 (C)]) \longrightarrow I/I^n(C).$$
We claim that for all $n>1$, this map is a weak equivalence. We induct on $n$.
 
For $n=2$, the left-hand side is equivalent to $I/I^2 (C)$, as clearly is the right-hand 
side. Moreover the map $I/I^2 (\alpha)$ is an equivalence since by assumption the composite
$I/I^2 (C) \to U_{\a}(C)  \to I/I^2 (C)$ is an equivalence. 

Next assume that the claim is true for $n=k$, i.e that 
$I/I^k (\alpha):I/I^k(T_{\a} [I/I^2 (C)]) \to I/I^k(C)$ is an equivalence, 
and consider the following 
diagram. Observe that the diagram is commutative since all the vertical maps are induced 
by $\alpha$. 
\begin{equation}
\label{diag:fib}
\xymatrix{
I^k/I^{k+1}(T_{\a} [I/I^2 (C)]) \ar[r] \ar[d]
&
I/I^{k+1}(T_{\a} [I/I^2 (C)]) \ar[r] \ar[d] 
&
I/I^k(T_{\a} [I/I^2 (C)]) \ar[d]\\
I^k/I^{k+1}(C) \ar[r]
&
I/I^{k+1}(C) \ar[r]
&
I^k(C)
}
\end{equation}
The right vertical map is a weak equivalence by inductive hypothesis. The source and the 
target of the left vertical map are both equivalent to 
${\a}(k) \wedge_{h \Sigma_k} \bigwedge ^k I/I^2(C)$ by Proposition~\ref{prop:fiber}. Moreover,
the proof of Proposition~\ref{prop:fiber} implies that the diagram
$$\xymatrix{
{\a}(k) \wedge_{h \Sigma_k} \bigwedge ^k I/I^2(T_{\a} [I/I^2 (C)]) \ar[r]^-{\simeq} \ar[d]
&
I^k/I^{k+1}(T_{\a} [I/I^2 (C)]) \ar[d]\\
{\a}(k) \wedge_{h \Sigma_k} \bigwedge ^k I/I^2(C) \ar[r]^-{\simeq}
&
I/I^k(T_{\a} [I/I^2 (C)])
}$$
commutes up to a weak equivalence. Here the left hand side map is induced by $\alpha$. 
We already argued that $I/I^2$ applied to $\alpha$ is an equivalence, and hence 
the left vertical map is an equivalence. 
This implies that the right vertical map, or equivalently (given that it is the same map), 
the left vertical map of Diagram~\ref{diag:fib} is an equivalence.  
Hence the middle vertical map of Diagram~\ref{diag:fib} is also a weak equivalence, thus 
proving that for all $n>1$, 
$$I/I^n (\alpha):I/I^n(T_{\a} [I/I^2 (C)]) \longrightarrow I/I^n(C)$$
is a weak equivalence. 
Consequently,
$$\holim_n I/I^n(C) \simeq \holim_n I/I^n(T_{\a} [I/I^2 (C)]) \simeq 
\prod {\a}(n) \wedge_{\Sigma_n} [I/I^2 (C)]^{\wedge n}$$

\end{proof}

Observe that whenever in addition to the hypothesises of this theorem, $C$ and $T_{\a} C$ 
are complete. i.e. they are equivalent to the homotopy limit of their respective 
Goodwillie towers (for example if the  $\a$-algebras $C$ 
satisfies the connectivity assumptions of Proposition~\ref{prop:conn}), this theorem
implies that $C$ is a free object in the category $\a$-algebras, thus answering the 
earlier posed question.

\begin{remark}
Naturally, Proposition~\ref{prop:genfiber} allows us to prove the analogue of this theorem 
for general operads. We do state the general version of the theorem for completeness of the 
picture and reference purposes, but do not present the proof here. Instead we comment that the 
proof of Theorem~\ref{theor:algsplit} works here nearly verbatim, once we replace $I/I^n (C)$ by 
$P_{n-1} U_{\a} (C)$, $I^n/I^{n+1} (C)$ by $D_{n} U_{\a} (C)$, and use  
Proposition~\ref{prop:genfiber} whenever Proposition~\ref{prop:fiber} is used.
\end{remark}

\begin{theorem}
\label{theor:genalgsplit}
Let the operad $\a$ be such that the operadic multiplications 
$${\a}(n) \wedge {\a}(1) \wedge \cdots \wedge {\a}(1) \to {\a}(n)$$
are equivalences. Furthermore, suppose the $\a$-algebra $C$ is such that the natural 
derivative map $U_{\a} (C) \to D_1 U_{\a}(C)$ has a section in the category of $A$-modules. 
Then we have a weak equivalence 
$$\holim_n P_n U_{\a} (C) \simeq \prod {\a}(n) \wedge_{h\Sigma_n} [D_1 U_{\a} (C)]^{\wedge n}.$$ 
\end{theorem}

Combining the results of this section with the main theorem of~\cite{Krist} yields an immediate 
corollary on the structure of $co-H-objects$ in the category of algebras over a fixed operad 
$\a$. We begin by recalling a definition from~\cite{Krist}.

\begin{definition}
Let ${\mathcal C}$ be a pointed model category, and let $X$ be a cofibrant object in 
${\mathcal C}$. We say that $X$ is a $co-H-object$ of ${\mathcal C}$ if there exists a
comultiplication $\nabla: X \to X \vee X$ which is coassociative and counital up to homotopy,
\end{definition}

The obvious relationship of $co-H-objects$ with the classical notion of quasi-Hopf algebras 
is what prompted the choice of terminology in the above definition. 

\begin{theorem}[Bauer-McCarthy,~\cite{Krist}]
\label{theor:krist}
Let $F$ be a functor from ${\mathcal C}$ to the category of spectra and let $X$ be 
a $co-H-object$ of ${\mathcal C}$. Then rationally,
$$\holim_n P_n F(X) \simeq \prod_{n \geq 0} D_n F(X).$$  
\end{theorem}

In particular, applying this theorem to the forgetful functor 
$U_{\a}: \C_{\a} \to \M_A$ for a fixed operad $\a$ as above, we conclude that if 
$C$ is a  $co-H-object$ of $\C_{\a}$, then the natural derivative map 
$U_{\a}(C) \to D_1 U_{\a}(C)$ has a section. Thus, the following corollary is immediate 
from Theorem~\ref{theor:algsplit} and Theorem~\ref{theor:genalgsplit}.

\begin{corollary}
\label{cor:krist}
Let the operad $\a$ be such that the operadic multiplications 
$${\a}(n) \wedge {\a}(1) \wedge \cdots \wedge {\a}(1) \to {\a}(n)$$
are equivalences. Furthermore, suppose the $\a$-algebra $C$ is  a  $co-H-object$ of $\C_{\a}$.
Then rationally we have a weak equivalence 
$$\holim_n P_n U_{\a} (C) \simeq \prod {\a}(n) \wedge_{h\Sigma_n} [D_1 U_{\a} (C)]^{\wedge n}.$$ 
If in addition, $\a$ is primitively generated, then rationally
$$\holim_n P_n U_{\a} (C) \simeq \holim_n I/I^n (C) 
\simeq \prod {\a}(n) \wedge_{\Sigma_n} [I/I^2 (C)]^{\wedge n}.$$
\end{corollary} 

Furthermore, in~\cite{Krist} the authors extend the above Theorem~\ref{theor:krist}
to the integral setting as follows.
\begin{theorem}
Let $F$ be a homotopy functor from ${\mathcal C}$ to the category of spectra and let $X$ be 
a cocommutative $co-H-object$ of ${\mathcal C}$. Then
$$\holim_n P_n F(X) \simeq \prod_{n \geq 0} D_n F(X)$$
whenever the Tate cohomology 
$$
{\Tate} ^n (F;X) \stackrel{def}{=} \cofiber[(D_1 ^{(n)} cr_n F(X))_{h \Sigma_n} \to 
(D_1 ^{(n)} cr_n F(X))^{h \Sigma_n}]
$$
vanishes for all $n$. 
\end{theorem}

Just as above, we combine this result with Theorem~\ref{theor:algsplit} and 
Theorem~\ref{theor:genalgsplit} to produce a corollary analogous to 
Corollary~\ref{cor:krist}, which we state here for completeness and reference purposes.
\begin{corollary}
Let the operad $\a$ be such that the operadic multiplications 
$${\a}(n) \wedge {\a}(1) \wedge \cdots \wedge {\a}(1) \to {\a}(n)$$
are equivalences. Furthermore, suppose the $\a$-algebra $C$ is  a  
cocommutative $co-H-object$ of $\C_{\a}$.
Then whenever the Tate cohomology $\Tate^n(U_{\a};C)$ vanishes for all $n$,  
we have a weak equivalence 
$$\holim_n P_n U_{\a} (C) \simeq \prod {\a}(n) \wedge_{h\Sigma_n} [D_1 U_{\a} (C)]^{\wedge n}.$$ 
If in addition, $\a$ is primitively generated, then 
$$\holim_n P_n U_{\a} (C) \simeq \holim_n I/I^n (C) 
\simeq \prod {\a}(n) \wedge_{\Sigma_n} [I/I^2 (C)]^{\wedge n}.$$
\end{corollary}

\subsection{Some Classical Theorems}

As an example we consider the case when the operad $\a$ is an $E_\infty$ operad, 
in other words the category of $\a$-algebras is the category of commutative $A$-algebras.  
In this case, whenever the commutative $A$-algebra $C$ is such that the conditions of 
the theorem are satisfied we get the identity
$$ \holim_n I/I^n(C) \simeq \P_A [I/I^2 (C)],$$
where $\P$ is the symmetric algebra cotriple, i.e. 
$\P(M) \cong \bigoplus M^{\wedge n}/\Sigma_n$, with the smash products taken over $A$.
Moreover, if $C$ is the
Topological Hochschild homology of the commutative $A$-algebra $A$, i.e. 
by a theorem of McClure, Schw\"{a}nzl and Vogt (\cite{Vogt}), $C \simeq A \otimes S^1$,
then $I/I^2(C)$ is equivalent to the suspension of Andr\'e-Quillen homology, in other 
words $I/I^2(C) \simeq \Sigma TAQ(A)$ (see e.g.~\cite{Minas}). Hence, under the 
connectivity conditions of Proposition~\ref{prop:conn}, we get
$$THH(A) \simeq \P_A[\Sigma TAQ(A)].$$
The next obvious question is `For which algebras $A$ are the conditions of the theorem 
satisfied for $THH(A)$?' This has been discussed in detail in~\cite{McMin}. Here we 
present the answer and refer to ~\cite{McMin} for proofs. We start by recalling 
a couple of definitions from there.

\begin{definition}
The map of commutative $S$-algebras $C \to D$ is \'etale if $TAQ(D|C)$ is contractible.

We say that $\{ C \to C_\alpha \}_{\alpha \in \I}$ is an \'etale covering of $C$ 
if

1. each map  $ C \to C_\alpha $ is \'etale, and 

2. for each pair of $C$-modules $M \to N$ such that 
$M \wedge C_\alpha \to N \wedge C_\alpha$ has a section for all $\alpha$, 
the map $M \to N$ itself has a section $N \to M$ with $N \to M \to N$ equivalent 
to the identity on $M$. 
\end{definition}

\begin{definition}
The map of algebras $f:R \to C$ is smooth if there is an \'etale covering 
$\{ C \to C_\alpha \}_{\alpha \in \I}$ of $C$ such that for each $\alpha$ there is 
a factorization
$$ R \longrightarrow \P_R X \stackrel{\phi}{\longrightarrow} C_\alpha, $$
where $X$ is a cell $R$-module and $\P_R X$ is the free commutative $R$-algebra 
generated by $X$, with $\phi$ \'etale. 
\end{definition}

The following theorem which answers our earlier posed question, has its analogue
in discrete algebra, where it is commonly referred to as 
Hochschild-Kostant-Rosenberg (HKR) theorem. For complete details in discrete case
we refer to~\cite{Loday}.

\begin{theorem}
(HKR theorem,~\cite{McMin})
Let $f:R \to A$ be smooth in the category of connective 
$S$-algebras. Then the natural (derivative) map 
$THH(A|R) \to \Sigma TAQ (A|R)$ has a section in the category of $A$-modules 
which induces an equivalence of $A$-algebras:
$$\P_A \Sigma TAQ (A|R) \stackrel{\simeq}{\longrightarrow} THH(A|R).$$
\end{theorem} 

Next we illustrate how the theorem of Leray on structure of commutative quasi Hopf algebras 
is a special case of Theorem~\ref{theor:algsplit}. First we recall some basic 
constructions from the theory of Hopf algebras. The extensively used paper~\cite{Milnor}
of J.Milnor and  J.Moore is our main reference. 

Let $K$ be a commutative ring, and $A$ and augmented $K$-algebra. We denote by $I(A)$ the 
augmentation ideal of $A$, i.e. $I(A)\stackrel{def}{=} \ker [A \to K]$. 
\begin{definition} 
If $A$ is an augmented algebra over $K$, let $Q(A)= K \otimes_A I(A)$. The elements of the 
$K$-module $Q(A)$ are called the indecomposables elements of $A$. Further, note that there 
is a natural exact sequence 
$$
I(A) \otimes I(A) \to I(A) \to Q(A) \to 0.
$$                    
\end{definition}
The following is listed as Theorem 7.5 in~\cite{Milnor}.
\begin{theorem}
(Leray)
If $A$ is a connected commutative quasi Hopf algebra over the field $K$ of characteristic 
zero and $X =Q(A)$, then if $f:X \to I(A)$ is a morphism of graded vector spaces such that 
the composition $X \stackrel{f}{\to} I(A) \to X$ is the identity morphism of $X$, then there 
is an isomorphism of algebras $\P X \to A$ induced by $f$, where $\P X$ is the polynomial 
algebra on $X$.   
\end{theorem}

Connectedness of $A$ is assumed to guarantee that the augmentation filtration on 
$A$ is complete. 

Observe that a connected commutative quasi Hopf algebra is in particular a commutative algebra.
Thus, by Theorem~\ref{theor:algsplit} (or more precisely the spacial case of 
${\a}={\bf e_\infty}$), the algebra $A$ is equivalent to the symmetric algebra on $I/I^2(A)$
provided the derivative map $A \to I/I^2(A)$ has a section. However, the existence 
of such a section is precisely the hypothesis of the Leray Theorem once we recall that the 
module of indecomposables $Q(A)$ is  $I/I^2(A)$. Thus the Leray theorem is simply 
a special case of our Theorem~\ref{theor:algsplit}.

\section{Operads Induced by Triples}
\label{sect:tripop}

Earlier, we recalled a construction that associated a triple to the given operad. 
In particular, it was easily seen that the algebras over the operad are same as the 
algebras over the triple. 

In this section we introduce a construction that produces operads out of triples in the 
category of spectra. 

It is worth noting that Andrew Mauer-Oats is currently working independently on a similar 
construction for functors of spaces, which presents challenges different (and most likely more 
complex) than the ones we encountered here while working with functors of spectra. 

The notions of cross effects $cr_n F$ of a functor 
$F: {\mathcal B} \to {\mathcal A}$ from a basepointed category $ {\mathcal B}$ to an additive
category $ {\mathcal A}$, and of the derivative (or linearization) functor $D_1 F$ are 
the main ingredients of our constructions. See Section~\ref{sec:derdif} for a brief discussion
on both. 

Recall that roughly speaking,  $D_1 F$ is the linear approximation of $F$. It comes equipped 
with a natural (derivative) map $F \to D_1 F$ which is an equivalence if $F$ is already of 
degree 1 (or linear). In particular, $D_1 D_1 F$ is equivalent to $D_1 F$. 

The following Chain Rule lemma is a restatement in our set up of Lemma 5.7 of~\cite{Randy}.
\begin{lemma}[Chain Rule]
\label{lem:chain}
Let $F$ and $G$ be endofunctors of the category of $S$-modules $\M _S$. Then 
$$D_1 (F \circ G) \simeq D_1 F \circ D_1 G.$$
\end{lemma} 

We recall our notation convention. For a 
functor $F$ of $n$ variables, we denote by $D_1 ^i F(X_1, \cdots , X_n)$ the derivative of $F$ 
obtained by holding all but $i'th$ variable constant. Also, we denote 
by $D_1 ^{(n)}F$ the multilinearized functor $D_1 ^n \cdots D_1 ^2 D_1 ^1 F$.

Let $F:\M_S \to M_S$ be a triple in the category of $S$-modules. In other words 
we have a multiplication map $\mu: F \circ F \Rightarrow F$ and a unit map 
$\eta: id_{\M_S} \Rightarrow F$, that satisfy the usual associativity and unit 
diagrams. Define
$${\a}_F (n) = D_1 ^{(n)} cr_n F(S).$$  
Our next objective is to demonstrate that the collection $\{{\a}_F (n) \}$ forms an operad.
Naturally, we need to produce the operadic multiplication maps. 

First note that since $id_{\M _S} \stackrel{\simeq}{\to} D_1 id_{\M _S}$,  applying the 
derivative functor $D_1$ to the natural transformation $\eta$
produces the unit map $S \to D_1 cr_1 F(S)$. In what follows, for simplicity, we are going to assume 
that $F$ is reduced, i.e $cr_1 F =F$.

Now fix an object $X\in \M_S$ and
consider the triple multiplication map $F \circ F (\vee_n X) \to F(\vee_n X)$. We apply 
the functor $D_1 ^1$ (i.e. the derivative is taken with respect to the first summand, while 
the others are held constant), and use the Chain Rule Lemma~\ref{lem:chain} to get 
$$D_1 ^1 (F \circ F (\vee_n X)) \simeq D_1 F \circ D_1 ^1 F(\vee_n X) \to D_1 ^1 F(\vee_n X).$$  
Iterate this construction 
by successively applying $D_1 ^2, D_1 ^3, \cdots , D_1 ^n$, to get 
\begin{equation}
\label{eq:opmult1}
D_1 F \circ D_1 ^{(n)} F(\vee_n X) \to D_1 ^{(n)} F(\vee_n X).
\end{equation}
Recall that $D_1 F$ is a reduced linear functor. Hence, for each $Y \in \M_S$, 
$D_1 F(Y) = D_1 F(S \wedge Y) \cong D_1 F(S) \wedge Y$. Consequently, noting that 
$ D_1 ^{(n)} F(\vee_n X) \simeq D_1 ^{(n)} cr_n F(X)$, we can rewrite the above 
Equation~\ref{eq:opmult1}:
\begin{equation}
\label{eq:opmult2}
D_1 F(S) \wedge D_1 ^{(n)} cr_n F(X) \to  D_1 ^{(n)} cr_n F(X). 
\end{equation}
Thus taking $X=S$ in this equation, we get the multiplication maps 
${\a}_F (1) \wedge {\a}_F (n) \to {\a}_F (n)$.
 
Before producing the general multiplication map, it is beneficial (for better clarity) 
to understand the maps ${\a}_F (2) \wedge {\a}_F (n) \wedge {\a}_F (m) \to {\a}_F (n+m)$
first. To that end, consider the composite map
\begin{equation}
\label{eq:opmult3}
F \circ [F(X_1 \vee \cdots \vee X_n) \vee F(X_{n+1} \vee \cdots \vee X_{n+m})] \to 
F \circ  F(X_1 \vee \cdots \vee X_{n+m}) \stackrel{\mu}{\to} F(X_1 \vee \cdots \vee X_{n+m}),
\end{equation}
where the first map is induced by the obvious inclusions 
$F(X_1 \vee \cdots \vee X_n) \to F(X_1 \vee \cdots \vee X_{n+m})$ and 
$F(X_{n+1} \vee \cdots \vee X_{n+m}) \to F(X_1 \vee \cdots \vee X_{n+m})$.
Applying $D_1 ^1, D_1 ^2, \cdots , D_1 ^{n}$ successively, using the chain 
rule (Lemma~\ref{lem:chain}) again, and setting $X_1= \cdots =X_{n+m} =X$, we get 
\begin{equation}
\label{eq:opmult4}
D_1 ^1 F [D_1 ^{(n)} F(\vee _n X) \vee F(\vee _m X)] 
\to D_1 ^{(n)} F(\vee _{n+m} X),
\end{equation}
where on the right the derivatives are taken with respect to the first $n$ copies of $X$.
Here we are treating the leftmost $F$ of Equation~\ref{eq:opmult3}
as a functor of two variable in the obvious way, consequently, the superscript $1$ in 
$D_1 ^1 F$ on the left of the Equation~\ref{eq:opmult4} is a reflection of the fact 
that all $n$ derivatives are taken over the first component $F(X_1 \vee \cdots \vee X_n)$.
Further applying $D_1 ^{n+1}, \cdots , D_1 ^{n+m}$, using the 
fact that 
$ D_1 ^{(l)} F(X_1 \vee \cdots \vee X_l) \cong D_1 ^{(l)} cr_l F(X_1, \cdots, X_l)$ 
for all $l>1$, and continuing to set $X_1= \cdots =X_{n+m} =X$, we get a map 
\begin{equation}
\label{eq:opmult5}
D_1 ^{(2)} cr_2 F [D_1 ^{(n)} cr_n F(X), D_1 ^{(m)} cr_m F(X)] \to 
D_1 ^{(n+m)} cr_{n+m} F(X).
\end{equation}
Here we employed the fact that the last $m$ derivatives are taken over the second component 
$F(X_{n+1} \vee \cdots \vee X_{n+m})$ of the left hand side of Equation~\ref{eq:opmult3}.

Recall that $D_1 ^{(2)} cr_2 F$ is linear in each of its variables. Hence, for $S$-modules 
$M$ and $N$, we have an identity 
$$D_1 ^{(2)} cr_2 F(M,N) \cong D_1 ^{(2)} cr_2 F(S,S) \wedge M \wedge N = 
D_1 ^{(2)} cr_2 F(S) \wedge M \wedge N.$$ 
Thus, Equation~\ref{eq:opmult5} can be rewritten as 
\begin{equation}
\label{eq:opmult6}
D_1 ^{(2)} cr_2 F(S) \wedge D_1 ^{(n)} cr_n F(X)\wedge D_1 ^{(m)} cr_m F(X) \to 
D_1 ^{(n+m)} cr_{n+m} F(X).
\end{equation}
Specializing to $X=S$, we get the desired multiplication map 
${\a}_F (2) \wedge {\a}_F (n) \wedge {\a}_F (m) \to {\a}_F (n+m)$.

The general case of the map 
${\a}_F(k) \wedge {\a}_F (j_1) \wedge \cdots \wedge {\a}_F (j_k) \to
{\a}_F (j_1 +\cdots +j_k)$
is analogous in spirit. Similar to the previous case, we employ the map 
\begin{eqnarray*}
F \circ [F(X_1 \vee \cdots \vee X_{j_1}) \vee \cdots \vee  
F(X_{j_1 + \cdots + j_{k-1} +1} \vee \cdots \vee X_{j_1 + \cdots + j_k})] 
& \to & 
F\circ F(X_1 \vee \cdots \vee X_{j_1 + \cdots + j_k})\\
& \to & 
F(X_1 \vee \cdots \vee X_{j_1 + \cdots + j_k})
\end{eqnarray*}
induced by evident inclusions 
$X_{j_1 + \cdots + j_{l-1} +1} \vee \cdots \vee X_{j_1 + \cdots + j_l} \to 
X_1 \vee \cdots \vee X_{j_1 + \cdots + j_k}$ (with $l \leq k$). Applying the sequence of 
$D_1 ^i$'s and setting $X_j =X$ for all $j$, we get the analogue of Equation~\ref{eq:opmult5}:
\begin{equation}
\label{eq:opmult7}
D_1 ^{(k)} cr_k F[D_1 ^{(j_1)} cr_{j_1} F(X), \cdots, D_1 ^{(j_k)} cr_{j_k} F(X)] \to 
D_1 ^{(j_1 + \cdots + j_k)} cr_{j_1 + \cdots + j_k} F(X).
\end{equation}
Finally, noting that $D_1 ^{(k)} cr_k F$ is linear in each of its $k$ variables,
we rewrite Equation~\ref{eq:opmult7} as
\begin{equation}
\label{eq:opmult8}
D_1 ^{(k)} cr_k F (S) \wedge D_1 ^{(j_1)} cr_{j_1} F(X) \wedge \cdots \wedge  
D_1 ^{(j_k)} cr_{j_k} F(X) \to D_1 ^{(j_1 + \cdots + j_k)} cr_{j_1 + \cdots + j_k} F(X),
\end{equation}
to produce the general case of the multiplication map by setting $X=S$. 

\begin{theorem}
\label{th:tripop}
With the above defined multiplication, the collection 
$\{ {\a}_F (n) = D_1 ^{(n)} cr_n F(S) \}_{n=0} ^{\infty}$ forms an operad.
\end{theorem}

\begin{proof}
We need to show that the associativity (Definition~\ref{def:operad}, part (a)), unital (part (b)) 
and equivariance (part (c)) diagrams commute for our operadic multiplication. We will use 
the same indexing as in Definition~\ref{def:operad}.

To prove the associativity consider the following object:
\begin{equation}
\label{object1}
F \circ [F \circ (F(\vee_{i_{1}} X) \vee \cdots \vee F(\vee_{i_{j_1}} X)) \vee \cdots \vee 
F \circ (F(\vee_{i_{g_{k-1}+1}} X) \vee \cdots \vee F(\vee_{i_{g_k}} X))].
\end{equation}
There are two maps from this object to 
$F(\vee_{i_{1}+\cdots +i_{j_1}+ \cdots +i_{g_{k-1}+1}+ \cdots +i_{g_{k}}} X) =
F(\vee_{i_{1}+\cdots + i_{g_{k}}} X)$.

{\bf Map 1.} There is an evident map from Object~\ref{object1} to 
\begin{equation}
\label{object2}
F \circ [F \circ F(\vee_{i_{1}+\cdots +i_{j_1}} X) \vee \cdots \vee 
F \circ F(\vee_{i_{g_{k-1}+1}+\cdots + i_{g_{k}}} X)],
\end{equation}
given by inclusions of $F(\vee_{i_{g_{l-1}+1}} X) \vee \cdots \vee F(\vee_{i_{g_l}} X)$
into $F(\vee_{i_{g_{l-1}+1}+\cdots + i_{g_{l}}} X)$ for $l=1, \cdots, k $. Furthermore,
noting that there is an evident inclusion map from 
$F \circ F(\vee_{i_{1}+\cdots +i_{j_1}} X) \vee \cdots \vee 
F \circ F(\vee_{i_{g_{k-1}+1}+\cdots + i_{g_{k}}} X)$ to 
$F \circ F (\vee_{i_{1}+\cdots +i_{j_1}+ \cdots +i_{g_{k-1}+1}+ \cdots +i_{g_{k}}} X)$,
we get a map from Object~\ref{object2} to 
\begin{equation}
\label{object3}
F\circ [F \circ F(\vee_{i_{1}+\cdots +i_{j_1}+ \cdots +i_{g_{k-1}+1}+ \cdots +i_{g_{k}}} X)]
= F\circ [F \circ F (\vee_{i_{1}+\cdots + i_{g_{k}}} X)]
\end{equation}
Finally, using the natural multiplication 
$F \circ (F \circ F) \Rightarrow F \circ F \Rightarrow F$, we get the desired map from 
Object~\ref{object1} to $F(\vee_{i_{1}+\cdots + i_{g_{k}}} X)$.

{\bf Map2.} Now we begin by considering the map from Object~\ref{object1} to 
\begin{equation}
\label{object4}
F \circ F \circ [F(\vee_{i_{1}} X) \vee \cdots \vee F(\vee_{i_{j_1}} X) \vee \cdots \vee 
F(\vee_{i_{g_{k-1}+1}} X) \vee \cdots \vee F(\vee_{i_{g_k}} X)],
\end{equation}
given by the inclusion of 
$F \circ (F(\vee_{i_{1}} X) \vee \cdots \vee F(\vee_{i_{j_1}} X)) \vee \cdots \vee 
F \circ (F(\vee_{i_{g_{k-1}+1}} X) \vee \cdots \vee F(\vee_{i_{g_k}} X))$ into 
$F \circ [F(\vee_{i_{1}} X) \vee \cdots \vee F(\vee_{i_{j_1}} X) \vee \cdots \vee 
F(\vee_{i_{g_{k-1}+1}} X) \vee \cdots \vee F(\vee_{i_{g_k}} X)]$. Moreover, using an analogous 
inclusion we get a map from Object~\ref{object4} to 
\begin{equation}
\label{object5}
(F\circ F) \circ [F(\vee_{i_{1}+\cdots +i_{j_1}+ \cdots +i_{g_{k-1}+1}+ \cdots +i_{g_{k}}} X)]
= (F\circ F) \circ [F (\vee_{i_{1}+\cdots + i_{g_{k}}} X)].
\end{equation}
Similar to the previous case, the natural multiplication 
$(F \circ F) \circ F \Rightarrow F \circ F \Rightarrow F$, produces the second map from 
Object~\ref{object1} to $F(\vee_{i_{1}+\cdots + i_{g_{k}}} X)$.

Of course, our interest in these maps is explained by the fact that applying the 
multilinearization functor $D_1 ^{(i)}$ to the target and source and setting $X=S$,
we get multiplication maps
\begin{equation*}
{\a}_F (k) \wedge (\bigwedge _{s=1} ^k {\a}_F (j_s)) \wedge 
(\bigwedge _{r=1} ^j {\a}_F (i_r)) \to {\a}_F (i),
\end{equation*}
where the indexing is still the same as in Definition~\ref{def:operad}. Moreover, note that 
$D_1 ^{(i)}$ applied to Map 1 is the composition of the smash product of multiplications
\begin{equation*} 
{\a}_F (j_l) \wedge {\a}_F (i_{g_{l-1}+1}) \wedge \cdots \wedge {\a}_F (i_{g_l}) \to 
{\a}_F (i_{g_{l-1}+1} + \cdots + i_{g_l})
\end{equation*}
for $l=1, \cdots k$, with the map
\begin{equation*}
{\a}_F (k) \wedge {\a}_F(i_1 + \cdots + i_{g_1}) \wedge \cdots \wedge  
{\a}_F (i_{g_{k-1}+1} + \cdots + i_{g_k}) \to 
{\a}_F (i_1 + \cdots + i_{g_1} +\cdots +i_{g_{k-1}+1} + \cdots + i_{g_k}) = {\a}_F (i).
\end{equation*}
Observe that this composite corresponds to the left vertical map followed by the 
bottom horizontal map in the diagram of part (a) of Definition~\ref{def:operad}.

Similarly, $D_1 ^{(i)}$ of Map 2 is the composition of the maps 
\begin{equation*} 
{\a}_F (k) \wedge {\a}_F(j_1) \wedge \cdots \wedge {\a}_F(j_k) 
\to {\a}_F(j_1 +\cdots + j_k) = {\a}_F(j)
\end{equation*}
and 
\begin{equation*}
{\a}_F(j) \wedge {\a}_F(i_{1}) \wedge \cdots \wedge {\a}_F(i_{g_1})  \wedge \cdots \wedge
{\a}_F(i_{g_{k-1}+1}) \wedge \cdots \wedge {\a}_F (i_{g_k}) \to 
{\a}_F (i_1 + \cdots + i_{g_k})= {\a}_F (i).
\end{equation*}
Consequently, the multilinearization of $D_1 ^{(i)}$ of Map 2 corresponds to the 
top horizontal arrow followed by the right vertical map in the diagram of 
part (a) of Definition~\ref{def:operad}. 

To conclude the proof of the associativity, we note that the only difference between 
Map 1 and Map 2 (i.e. before we apply the multilinearization functor) is the order in which 
the three copies of the triple $F$ are composed, and hence the associativity of the 
composition of the triple $F$ implies that associativity of the operad ${\a}_F$.

To see that the unit diagrams for the operadic multiplication (part (b) of 
Definition~\ref{def:operad}) commute, consider the commutative diagram
$$
\xymatrix{
F(\vee _k S) \ar[r]^{\simeq} \ar[d]
&
F(\vee _k S)\\
F \circ F(\vee _k S) \ar[ur]_{\mu} 
&
}
$$
where the left vertical map is given by $\eta \circ id$ or $id \circ \eta$. In other 
words we have two diagrams here, and both of them commute since $F$ is a unital triple. 
Observe that the commutative diagrams we get by applying $D_1 ^{(k)}$ to these diagrams, 
are precisely the unital diagrams of Definition~\ref{def:operad}, thus proving that our 
operadic multiplication is unital.  

To prove that the operadic multiplication is equivariant we need to show that the two 
diagrams of part (c) of Definition~\ref{def:operad} commute. Again we use the indexing
and notation introduced there. Note that the left diagram is simply the 
multilinearization of the diagram
$$
\xymatrix{
F \circ (F (\vee _{j_1} S) \vee \cdots \vee F (\vee _{j_k} S)) \ar[r] \ar[d]
& 
F \circ (F (\vee _{\sigma(j_1)} S) \vee \cdots \vee F (\vee _{\sigma(j_k)} S)) \ar[d]\\
F \circ F(\vee _{j_1} S \vee \cdots \vee \vee _{j_k} S) \ar[d]
&
F \circ F(\vee _{\sigma(j_1)} S \vee \cdots \vee \vee _{\sigma(j_k)} S) \ar[d]\\
F(\vee _{j_1} S \vee \cdots \vee \vee _{j_k} S) \ar[r]
&
F(\vee _{\sigma(j_1)} S \vee \cdots \vee \vee _{\sigma(j_k)} S),
}
$$
hence, it is enough to show that this diagram commutes. However, the commutativity 
of this diagram is immediate by functoriality of $F$ and the fact that 
$F \circ F \Rightarrow F$ is a natural transformation. 

It remains to prove that the right diagram of part (c) of Definition~\ref{def:operad} 
commutes. To ease the notation we will consider the special case $k=1$. The general 
case is identical to the case $k=1$. To keep track of the group action, we denote 
the sphere $S$ by $S_1, S_2, \cdots , S_j$, in other words, all $S_i$'s are equal to $S$ 
and the subscripting is to make the action of $\tau \in \Sigma_j$ transparent. Consider 
the diagram
$$
\xymatrix{
F \circ F(S_1 \vee \cdots \vee S_j) \ar[r] \ar[d]
&
F \circ F(S_{\tau(1)} \vee \cdots \vee S_{\tau(j)}) \ar[d]\\
F(S_1 \vee \cdots \vee S_j) \ar[r]
&
F(S_{\tau(1)} \vee \cdots \vee S_{\tau(j)}).
}
$$
This diagram is commutative since $F \circ F \Rightarrow F$ is a natural transformation. 
The right equivariance diagram of part (c) of Definition~\ref{def:operad} follows once 
we apply the multilinearization functor $D_1 ^{(j)}$ to this diagram.
\end{proof}

\begin{remark}
In the above construction, as well as in Theorem~\ref{th:tripop}, we considered triples 
in the category of $S$-modules $\M_S$. This was done simply as a matter of convenience,
as everything in this section can be repeated verbatim for triples in the category $\M_A$ 
of $A$-modules, where $A$ is a cofibrant commutative $S$-algebra. In fact, the construction 
can be performed in any symmetric monoidal category that falls within the framework 
of~\cite{Randy}, in other words, all symmetric monoidal categories where the identity functor 
is linear and for which the linearization 
functor $D_1$ satisfies the properties employed here, i.e. the chain rule property, the 
idempotency, and the correspondence of compositions to smash products (or whatever the 
symmetric monoidal operation happens to be). 

The chain complexes $Ch(K)$ over a commutative ring $K$ is an example of such a category. 
\end{remark}

\begin{example}
Here we present a few important examples in the category $Ch(K)$ of chain complexes 
over a field $K$ of characteristic $0$. 

{\it Commutative algebras}. The symmetric algebra triple 
$S(X)=\bigoplus_n X^{\otimes n} /{\Sigma_n}$ produces the commutative algebras in the 
category $Ch(K)$. By~\cite{Randy} or~\cite{Minas2}, 
$D_1 ^{(n)}cr_n S(X)=X^{\otimes n} $. Hence, the operad ${\a}_S$ is such that 
${\a}_S(n) = K$. In other words, the operad induced by $S(X)$ is the $e_\infty$ operad. 
In this and following examples we use the notation of~\cite{Getz} when referring to 
specific operads. See~\cite{Getz} for details on these operads.    

{\it Associative algebras}. The associative algebras are algebras over the tensor algebra
triple $T(X)=\bigoplus_n X^{\otimes n}$. By computations performed in Section 6 
of~\cite{Minas2}, $D_1 ^{(n)}cr_n T(K)= K[\Sigma_n]$, implying that the associated 
operad ${\a}_T$ is the operad $e_1$. 

{\it Lie algebras}.
For the free Lie algebra triple $L$, by Section 7 of~\cite{Minas2}, the multilinearization 
of the cross effects $D_1 ^{(n)} cr_n L(K)$ is equivalent to the $n-1$-dimensional 
{\it Lie} representation of the symmetric group $\Sigma_n$. Consequently, the induced 
operad ${\a}_L$ is the {\it Lie} operad. 
 
{\it n-Poisson algebras, $n \geq 2$.} The free $n$-Poisson algebras are produced by the triple 
$$P_n (X) = S (\Sigma ^{1-n} L (\Sigma^{n-1} X)),$$ 
where $S$ and $L$ are the symmetric
and Lie algebra triples respectively, and $\Sigma$ is the suspension (shift) functor
in the category of chain complexes. By Lemma 10.2 of~\cite{Minas2}, 
$$D_1 ^{(k)} cr_k P_n(K) \cong \bigotimes _{j=1} ^{k-1}H^{\ast} (\bigvee_j S^{n-1}),$$
where on the right hand side we have the tensor product of cohomologies (with coefficients 
in K) of the wedge of $j$ copies of the $n-1$ sphere $S^{n-1}$. Hence by a result of 
F.Cohen (Lemma 6.2 of~\cite{Cohen}), the operad ${\a}_{P_n}$ for $n \geq 2$, is equivalent 
to the homology of little $n$-cubes operad $e_n$, which in turn is isomorphic to the 
$n$-{\it Poisson} operad $p_n$. 

\end{example}

\begin{remark}
Another simple observation is that this mechanism of producing operads out of triples 
is in fact a functor from the category of triples ${\mathcal T}$ to the category 
of operads ${\mathcal O}$. Indeed, if $\tau : F \to G$ is a natural transformation 
which is a morphism of triples, then it induces the necessary maps 
$D_1 ^{(n)} cr_n F(S) \to D_1 ^{(n)} cr_n G(S)$, which are compatible with operadic 
multiplications, since $\tau$ respects the triple multiplications of $F$ and $G$. 

Moreover, note that if ${\a}$ is an operad and $T_{\a}$ is the triple associated with it, 
then 
\begin{equation}
\label{eq:optripop}
D_1 ^{(n)} cr_n T_{\a}(S) = 
D_1 ^{(n)} [cr_n \bigoplus _{k=0} ^{\infty} {\a} (k) \wedge _{h\Sigma_k} S^{\wedge k}] \cong 
{\a}(n).
\end{equation}
To see this, recall that 
$D_1 ^{(n)} cr_n T_{\a}(X) \cong D_1 ^{(n)} T_{\a}(X \vee \cdots \vee X)$, where there 
are $n$ copies of $X$ on the right hand side. To keep better track, we will denote these 
$n$ copies $X_1, \cdots , X_n$. Consider the summand 
${\a} (k) \wedge _{h\Sigma_k} (X_1 \vee \cdots \vee X_n )^{\wedge k}$ for a fixed $k$. Expanding 
this term, we get a new sequence of summands each of which has precisely $k$ factors. Hence, 
if $k<n$ then from each summand at least one $X_i$ is missing, i.e. that summand is constant
with respect to $X_i$. Thus $D_i$ applied to that summand is contractible. Consequently,
$D_1 ^{(n)} {\a} (k) \wedge _{h\Sigma_k} (X_1 \vee \cdots \vee X_n )^{\wedge k}$ vanishes for 
all $k<n$. Now if $k>n$, then each of the summands has at least one $X_j$ in degree 2 or higher.
In other words, that summand is at least 2-multireduced (see Definition 3.1 of~\cite{Randy}) in 
that variable. Hence, by Proposition 3.2 of~\cite{Randy},    
$D_1 ^{(n)} {\a} (k) \wedge _{h\Sigma_k} (X_1 \vee \cdots \vee X_n )^{\wedge k}$ vanishes for
all $k>n$ as well. Finally,  consider the case $k=n$. By the same reasoning, all the summands 
with at least one $X_i$ missing, vanish. Hence, there are precisely $n!$ surviving terms - one 
for each way to form a string of length $n$ out of $X_1, \cdots , X_n$ without repeating any 
of $X_i$'s. Equation~\ref{eq:optripop} follows.  

In other words, if (by abuse of notation) we denote by ${\mathcal O}$ the category of operads 
on the full subcategory of cofibrant objects in $\M_S$ (as opposed on all of $\M_S$), 
we are allowed to conclude that the composite of functors 
${\mathcal O} \to {\mathcal T} \to {\mathcal O}$ is equivalent to the identity functor, as it 
is straightforward to check that it is the identity map on morphisms as well. 
The same, of course, is not the case with the composite 
$\alpha: {\mathcal T} \to {\mathcal O} \to {\mathcal T}$, as it is a well known fact 
that not every triple arises from an operad. In fact, the question of description of 
triples which are fixed under $\alpha$ is a splitting question, thus explaining our 
interest in above constructions in the context of this work. Indeed, if 
$F \in {\mathcal T}$ is such that $\alpha (F) \cong F$, then $F$ is the triple associated 
to some operad ${\a}$, and hence $F (X) \simeq \bigoplus _n {\a}(n) \wedge _{h\Sigma_n} X^{\wedge n}$.
\end{remark}

We discuss this next.

\section{Triples and Splitting}
\label{sect:splittrip}

It is the intent of this section to develop splitting criteria for triples 
$F:\M_S \to M_S$ in the category of $S$-modules. The work leading up to this point 
suggests two natural approaches. First, in Section~\ref{sect:splitfunc} we discussed 
the splitting of Goodwillie towers of functors landing in $\M_S$. Triples in 
$\M_S$, of course, are examples of such functors. Consequently, one could attempt to 
to specialize the results of Section~\ref{sect:splitfunc}. 

The second approach is suggested by splitting conditions for algebras over an operad
(Section~\ref{sect:splitalg}), since $F$ produces an operad ${\a}_F$, and thus one can 
hope that the results of Section~\ref{sect:splitalg} can be applied whenever $F(X)$ is 
an ${\a}_F$-algebra for some $X$ in $\M_S$. 

Here we explore questions arising as a result of pursuing these two directions.  
In this section we assume that $F$ is complete (unless 
explicitly specified otherwise), or that 
$F$ is equivalent to the homotopy limit of its Goodwillie tower. To avoid 
double subscripting, we denote the triple 
associated with the operad ${\a}_F$ by $T_F$ instead of $T_{{\a}_F}$, 
which would be more consistent with the notation used throughout this work.      

First, following on Definition V.2.1 of~\cite{EKMM}, we present a new piece of terminology.
\begin{definition}
For an operad ${\a}$ (or a triple $T$) in the category $\M_S$ of of $S$-modules, 
we say that $C$ is an ${\a}$-ring spectrum (or a $T$-ring spectrum) if it is equipped 
with multiplication maps ${\a}(n) \wedge_{\Sigma_n} C^{\wedge n} \to C$ (or $TC \to C$),
such that the usual unital and associativity diagrams commute up to a weak equivalence.   
\end{definition}

\begin{proposition}
\label{prop:trspnec}
Let ($F$, $\mu:F \circ F \Rightarrow F$, $\eta:Id \Rightarrow F$) be a triple in the category 
of $S$-modules, such that for all $X$ in $\M_S$, the Goodwillie tower of $F(X)$ splits.  
Then $F(X)$ is an ${\a}_F$-ring spectrum for all $X$.     
\end{proposition}

\begin{proof}
Since the tower of $F$ splits for all $X$ and $F$ is complete, each map 
$p_n :F(X) \to P_n F(X)$ has a 
section up to a weak equivalence, and each $P_n$ decomposes into the coproduct of layers,
i.e. $P_n \simeq \vee_{i=1} ^n D_i F(X)$. Consequently, for all $n$, we have a natural map
$D_n F(X) \to F(X)$ which, up to a weak equivalence, is a section to the composite
$F(X) \to P_n F(X) \to D_n F(X)$, where the second map is the evident projection. 
Moreover, by Proposition~\ref{prop:deriv} $D_n F$ is naturally equivalent 
to $D_1 ^{(n)} cr_n F_{h\Sigma_n}$. Hence we have natural maps 
$D_1 ^{(n)} cr_n F(X) \to F(X)$.

Recall that $D_1 ^{(n)} cr_n F$ is a short hand notation for 
$D_1 ^{(n)} cr_n F(X, \cdots , X)$, which is linear in each of its $n$ variables. 
Thus,
$$ D_1 ^{(n)} cr_n F(X) \simeq D_1 ^{(n)} cr_n F(S) \wedge X^{\wedge n}.$$
For more on this see Section~\ref{sect:splitfunc} or~\cite{Randy2}. Hence, we have 
natural maps $D_1 ^{(n)} cr_n F(S) \wedge_{\Sigma_n} X^{\wedge n} \to F(X)$
which by the universal property of coproducts induce a natural transformation (up to a 
weak equivalence)
\begin{equation}
\label{eq:tripmap}
\nu: T_F (X)=\bigoplus ^n D_1 ^{(n)} cr_n F(S) \wedge_{\Sigma_n} X^{\wedge n} 
\Rightarrow F(X).
\end{equation}
Note that both the source and the target of $\nu$ are triples, and the map itself 
respects the triple multiplications. In other words the diagrams

\begin{equation}
\label{diag:tripmap}
\xymatrix{
T_F \circ T_F (X) \ar[r] \ar[d]
&
T_F (X) \ar[d]\\
F \circ F (X) \ar[r]
&
F(X)
}
\end{equation}
commute for all $X$ (up to a weak equivalence). Indeed, recall that the triple 
multiplication of $T_F$ is, in essence, given by multilinearization (with respect 
to all variables) of the map 
\begin{eqnarray*}
F \circ [F(X_1 \vee \cdots \vee X_{j_1}) \vee \cdots \vee  
F(X_{j_1 + \cdots + j_{k-1} +1} \vee \cdots \vee X_{j_1 + \cdots + j_k})] 
& \to & 
F\circ F(X_1 \vee \cdots \vee X_{j_1 + \cdots + j_k})\\
& \to & 
F(X_1 \vee \cdots \vee X_{j_1 + \cdots + j_k}),
\end{eqnarray*}
which is the multiplication map of $F$.
See Section~\ref{sect:tripop} for details. 

Thus, we have produced a natural (up to a weak equivalence) map of 
triples $T_F \Rightarrow F$ which, in turn, gives rise to a multiplication
\begin{equation}
\label{eq:compat}
T_F \circ F(X) \stackrel{\nu_{F(X)}}{\longrightarrow} F \circ F(X) 
\stackrel{\mu}{\longrightarrow} F(X).
\end{equation}   
We claim that this map makes $F(X)$ into an ${\a}_F$-ring spectrum. To see this,
it remains to show that the associativity diagram
\begin{equation}
\label{diag:asstrip}
\xymatrix{
TT(F(X)) \ar[r] \ar[d]^=
&
T(F(F(X))) \ar[r]
&
F(F(F(X))) \ar[r]^{F(\mu)}
&
F(F(X)) \ar[d]\\
TT(F(X)) \ar[r]
&
T(F(X)) \ar[r]
&
F(F(X)) \ar[r]
&
F(X)
}
\end{equation} 
commutes (up to a weak equivalence). Note that replacing $X$ by $F(X)$ in the commutative 
Diagram~\ref{diag:tripmap}, we have that the bottom horizontal arrow is equivalent to 
$$
TT(F(X)) \to FF(F(X)) \stackrel{\mu_{F(X)}}{\to} F(F(X)) \to F(X). 
$$
We rewrite Diagram~\ref{diag:asstrip} as
\begin{equation*}
\xymatrix{
TT(F(X)) \ar[r] \ar[d]^=
&
T(F(F(X))) \ar[r]
&
F(F(F(X))) \ar[r]^{F(\mu)}
&
F(F(X)) \ar[d]\\
TT(F(X)) \ar[r]
&
FF(F(X)) \ar[r]^{\mu_{F(X)}} \ar[ru]^{\simeq}
&
F(F(X)) \ar[r]
&
F(X),
}
\end{equation*} 
where the commutativity of the left half is evident, and the right half commutes 
by associativity of the triple multiplication of $F$.

\end{proof}

Note that in the above proposition, $F(X)$ is equipped with two multiplications, 
namely $T_F \circ F(X) \to F(X)$ and $F \circ F(X) \to F(X)$, and they are compatible 
in the sense that the diagram 
\begin{equation}
\label{diag:compat}
\xymatrix{
T_F \circ F(X) \ar[r] \ar[d]
&
F(X) \ar[d]^=\\
F \circ F(X) \ar[r]
&
F(X)
}
\end{equation}
commutes (up to a weak equivalence).
This is a consequence of the way the ${\a}_F$-ring spectrum multiplication of $F(X)$ 
was defined; see Equation~\ref{eq:compat}.

Our immediate goal is to show that the converse of Proposition~\ref{prop:trspnec} 
also holds. First we make a few comments. 

Suppose we have a triple $F$ which happens to have a natural structure of 
an ${\a}_F$-ring spectrum via a multiplication map 
$m:T_F(F(X)) \to F(X)$, or equivalently, via a sequence of maps 
$D_1 ^{(n)} cr_n F(S) \wedge_{\Sigma_n} F(X)^{\wedge n} \to F(X)$. 
Recalling again that $D_1 ^{(n)} cr_n F$ is linear in each of its $n$ variables, this 
maps produce a new set of morphisms 
\begin{equation}
\label{eq:helpmult}
D_1 ^{(n)} cr_n F(F(X))_{\Sigma_n} \to F(X).
\end{equation}
Furthermore, applying the multilinearized functor $D_1 ^{(n)} cr_n F_{\Sigma_n}$ 
to the unit $\eta:Id \Rightarrow F$ of the triple $F$, we get morphisms
$D_1 ^{(n)} cr_n F(X)_{\Sigma_n} \to D_1 ^{(n)} cr_n F(F(X))_{\Sigma_n}$,
which combine with Equation~\ref{eq:helpmult} to produce maps 
$D_1 ^{(n)} cr_n F(X)_{\Sigma_n} \to F(X)$.
By universal property of coproducts, these induce a natural transformation 
$\nu:T_F(X) \to F(X)$. 

This allows us to introduce a definition.

\begin{definition}
For any triple $F$, which is naturally an ${\a}_F$-ring spectrum, we say that the 
two algebra structures on $F(X)$ (given by $\mu: F \circ F(X) \to F(X)$ and
$\nu: T_F \circ F(X) \to F(X)$) are compatible if Diagram~\ref{diag:compat} commutes.  
\end{definition}

Now we are ready to state the converse of Proposition~\ref{prop:trspnec}. 

\begin{proposition}
Let ($F$, $\mu$, $\eta$) be a triple in $\M_S$ such that $F(X)$ is naturally an ${\a}_F$-ring 
spectrum for all $X$. Moreover, suppose the two algebra structures of $F(X)$ are 
compatible.  
Then the Goodwillie tower of $F$ splits for all $X$. 
\end{proposition}

\begin{proof}
Since $F$ is an ${\a}_F$-ring spectrum, by the discussion following 
Proposition~\ref{prop:trspnec}, we have a sequence of natural maps 
\begin{equation}
\label{eq:tripsect}
D_n F(X) \simeq D_1 ^{(n)} cr_n F(X)_{h\Sigma_n} \to D_1 ^{(n)} cr_n F(F(X))_{h\Sigma_n} 
\to D_1 ^{(n)} cr_n F(F(X))_{\Sigma_n} \to F(X).
\end{equation}
By inducting on $n$, we show that these maps split the Goodwillie tower of $F(X)$.

We begin with the base case $n=1$. Equation~\ref{eq:tripsect} specializes to give us a 
map $D_1 F(X) \to F(X)$ , which we need to show is a section (up to a weak equivalence) 
to the derivative map $F(X) \to D_1 F(X)$. Consider the following diagram.
\begin{equation}
\label{diag:d1split}
\xymatrix{
D_1 F(X) \ar[r] \ar[rrd] 
&
D_1 F(F(X)) \ar[r] \ar[rd]
&
F(X) \ar[r]^{p_1}
&
D_1 F(X)\\
&&
D_1 F(D_1 F(X)) \ar[ru]
&
}
\end{equation}
where the top vertical map is the case $n=1$ of Morphism~\ref{eq:tripsect}, and
the three non-horizontal maps are defined as follows. 
$D_1 F(X) \to D_1 F(D_1 F(X))$ is obtained by applying $D_1$ to the 
$F(\eta):F(X) \to F \circ F(X)$, and using the Chain Rule Lemma (see Lemma~\ref{lem:chain})
on the right hand side. The map $D_1 F(F(X)) \to D_1 F(D_1 F(X))$ is the functor 
$D_1 F$ applied to the derivative $F(X) \to D_1 F(X)$, and finally, 
$D_1 F(D_1 F(X)) \to D_1 F(X)$ is the derivative of the multiplication map 
$\mu : F \circ F \to F$. 

Note that our objective is to show that the top horizontal composite of 
Diagram~\ref{diag:d1split} is equivalent to the identity, and showing the commutativity 
of the diagram would accomplish that, since the lower composite 
$D_1 F(X) \to D_1 F(D_1 F(X)) \to D_1 F(X)$ is simply $D_1$ applied to 
$F(X) \stackrel{F(\eta)}{\longrightarrow} F \circ F \stackrel{\mu}{\longrightarrow} F$,
which is the identity map by the unit diagram for the triple $F$. Thus, if 
Diagram~\ref{diag:d1split} is commutative (up to a weak equivalence), the top 
horizontal composite is equivalent to the identity.

The commutativity of the left half of the diagram is evident from the description 
of the maps involved. To see that the right half of Diagram~\ref{diag:d1split} 
commutes consider
\begin{equation}
\label{diag:spright}
\xymatrix{
D_1 F(F(X)) \ar[r] 
&
T_F(F(X)) \ar[r]^-{deriv} \ar[d]
&
D_1 (T_F \circ F)(X) \ar[r]^-{\simeq} \ar[d]
&
D_1 (F\circ F)(X) \simeq D_1 F(D_1 F(X)) \ar[d]^{D_1 (\mu)}\\
&
F(X) \ar[r]^{deriv}
&
D_1 F(X) \ar[r]^=
&
D_1 F(X),
}
\end{equation}
where the left-most horizontal arrow is the inclusion of the degree 1 component 
$D_1 F(F(X))$ of $T_F \circ F$ into $T_F \circ F$. Here one should be careful not to 
confuse $D_1 F(F(X))$ with $D_1 (F\circ F)(X)$, as the former is $D_1 F$ evaluated at 
$F(X)$, while the latter is the derivative of the functor $F\circ F$. Observe that 
the left square of Diagram~\ref{diag:spright} commutes since the right vertical arrow is 
the derivative of the left one. The right square commutes by compatibility of the two algebra 
structures on $F(X)$ - note that it is simply $D_1$ applied to the compatibility 
diagram. The map  $D_1 (T_F \circ F)(X) \to D_1 (F\circ F)(X)$ is an equivalence 
since the derivatives $D_1$ of $T_F$ and $F$ are equivalent (recall the definition of 
$T_F$), and by compatibility of the algebra structures.   

Consequently, the composite of all top horizontal arrows with the right-most vertical map 
is equivalent to the composition of the two right-most slanted maps 
$D_1 F(F(X)) \to D_1 F(D_1 F(X)) \to D_1 F(X)$ of Diagram~\ref{diag:d1split}.

Suppose for all $n \leq k$ we have that $P_n F(X)$ is equivalent to the coproduct of 
layers $D_1 F(X) \vee \cdots \vee D_n F(X)$. Consider the diagram 
$$
\xymatrix{
D_{k+1} F(X) \simeq D_1 ^{(k+1)} cr_{k+1} F(X)_{h\Sigma_{k+1}} \ar[r]^-{d_{k+1}}
&
P_{k+1} F(X) \ar[r]^{q_{k+1}} \ar[d]^{q_{k+1}}
&
P_k F(X) \ar[d]^{q_{k}}\\
D_{k} F(X) \simeq D_1 ^{(k)} cr_{k} F(X)_{h\Sigma_{k}} \ar[r]^-{d_k} \ar[ur]
&
P_{k} F(X) \ar[r]^{q_{k}}
&
P_{k-1} F(X). \ar@/^1pc/[l]
}
$$
To see that $P_{k+1} F(X)$ splits we simply need to show that the composite 
\begin{equation}
\label{eq:eq84}
P_k F(X) \simeq P_{k-1} F(X) \vee D_1 ^{(k)} cr_{k} F(X)_{h\Sigma_{k}} \to 
F(X) \stackrel{p_{k+1}}{\to} P_{k+1} F(X) \stackrel{q_{k+1}}{\to} P_k F(X)
\end{equation}
is equivalent to the identity map. Here the first map on the component $P_{k-1} F(X)$ 
exists because by inductive hypothesis, $P_{k-1} F(X)$ is equivalent to the coproduct 
of layers $D_1 \vee \cdots \vee D_{k-1}$, and the maps on layers are defined via 
Equation~\ref{eq:tripsect}. As in Theorem~\ref{theor:splitfunc},
it is enough to prove that  Map~\ref{eq:eq84} is equivalent 
to the identity map on the component $D_1 ^{(k)} cr_{k} F(X)_{h\Sigma_{k}}$, because 
the identity on the component  $P_{k-1} F(X)$ follows by inductive hypothesis.
In other words, we need show that the composite
\begin{equation}
\label{eq:eq85}
D_1 ^{(k)} cr_{k} F(X)_{h\Sigma_{k}} \to F(X) \stackrel{p_{k}}{\to} P_k F(X), 
\end{equation}
which is the restriction of Map~\ref{eq:eq84} to $D_1 ^{(k)} cr_{k} F(X)_{h\Sigma_{k}}$,
is equivalent to the map  
$D_1 ^{(k)} cr_{k} F(X)_{h\Sigma_{k}} \stackrel{d_k}{\to} P_k F(X)$ of the Goodwillie
tower of $F(X)$. We form a diagram analogous to the one for the base case:
\begin{equation}
\label{diag:pksplit}
\xymatrix{
D_1 ^{(k)} cr_{k}F(X)_{h\Sigma_{k}} \ar[r] \ar[d] 
&
D_1 ^{(k)} cr_{k}F(F(X))_{h\Sigma_{k}} \ar[r] \ar[d]
&
F(X) \ar[r]^{p_k}
&
P_k F(X)\\
P_k F(X) \ar[r]^-{P_k F(\eta)}
&
P_k (F \circ F)(X) \ar[rru]_-{P_k(\mu)}
&&
}
\end{equation} 
where the left map in the top horizontal row is induced by the unit $\eta :Id \Rightarrow F$,
and the middle horizontal map is the ${\a}_F$-algebra map on $F(X)$. To describe the middle 
vertical map, observe that for all functors $F$ and $G$, there is a natural transformation 
$P_k F(G(X)) \to P_k (F \circ G)(X)$ induced by natural maps 
$cr_t F(G(X_1, \cdots, X_t)) \to cr_t (F \circ G)(X_1, \cdots, X_t)$. These are of course 
consequences of the maps $G(X_1) \vee \cdots \vee G(X_t) \to  G(X_1 \vee \cdots \vee X_t)$
resulting from the universal property of coproducts. Hence, we have a natural map 
$P_k F(F(X)) \to P_k (F \circ F)(X)$, and the middle vertical arrow is simply the composite 
$D_1 ^{(k)} cr_{k}F(F(X))_{h\Sigma_{k}} \to P_k F(F(X)) \to P_k (F \circ F)(X)$.

Observe that to complete the proof of the inductive step, it is enough to show that 
Diagram~\ref{diag:pksplit} commutes, since the composite of the left vertical arrow with the
bottom horizontal and the slanted arrows is precisely the map 
$D_1 ^{(k)} cr_{k} F(X)_{h\Sigma_{k}} \stackrel{d_k}{\to} P_k F(X)$. The left square of 
Diagram~\ref{diag:pksplit} commutes by naturality of Goodwillie towers, and the commutativity 
of the right triangle is a consequence of the compatibility of the two algebra structures on 
$F(X)$. More precisely, it follows from the fact that, analogous to the case $n=1$, the 
diagram 
\begin{equation*}
\xymatrix{
T_F(F(X)) \ar[r]^{p_k} \ar[d]
&
P_k (T_F \circ F)(X) \ar[r] \ar[d]
&
P_k (F\circ F)(X) \ar[d]^{P_k (\mu)}\\
F(X) \ar[r]^{p_k}
&
P_k F(X) \ar[r]^=
&
P_k F(X)
}
\end{equation*}
is commutative due to compatibility of the algebra structures. 
\end{proof}
 
At the beginning of this section we required that the triple $F$ be complete. Of 
course this was a rather critical assumption as it allowed us to construct maps 
into $F$ whenever we had compatible morphisms into Taylor polynomials $P_n F$ of 
$F$. In Section~\ref{sect:splitfunc} we observed that one way to remove this 
assumption is to replace $F$ by its completion $P_\infty F$ in every statement.
Naturally, we need to be more careful here as we have the additional structure of 
a triple on our functors, and in order to perform replacements we first need to verify 
that the completion $P_\infty F$ is also equipped with a triple multiplication.

\section{Cell {\a}-algebras}
\label{sect:last}

In this section, we prove the technical results which were utilized in 
Sections~\ref{sect:forget} and~\ref{sect:splitalg}. In essence, we lay the groundwork to overcome the 
problems created by the fact that $T_{\a} C$ is not a cofibrant $A$-module even when $C$
is. First we state the main results again.

\begin{theorem}
\label{theor:tech}
Let $C$ be a cofibrant $\a$-algebra and $\gamma:Y \to T_{\a}^n C$ a cell $A$-module approximation.

{\bf 1.}
Then the induced map 
$${\a}(i) \wedge_{\Sigma_i} Y^{\wedge i} \to {\a}(i) \wedge_{\Sigma_i} (T_{\a} ^n C)^{\wedge i} $$ 
is a weak equivalence for all $i>0$.

{\bf 2.}
The projection 
${\a}(i)  \wedge_{h\Sigma_i} (T_{\a} ^n C)^{\wedge i} \to 
 {\a}(i) \wedge_{\Sigma_i} (T_{\a} ^n C)^{\wedge i}$
is an equivalence for all $n$ and $i$. 
\end{theorem}

We follow the strategy employed by M.Basterra in~\cite{Maria} to treat the special case of the 
non-unital commutative algebra operad. She, in turn, relied heavily on~\cite{EKMM}. 
We recall a pair of definitions from~\cite{Maria}. 

\begin{definition}
Let $\F$ be the class of all $A$-modules of the form
$$A \wedge_S S \wedge_{\mathcal I} {\mathcal I}(i) \ltimes _G K,$$
where $K$ any $G$ spectrum indexed on the universe $U^i$ that has the homotopy type of a 
$G$-CW spectrum for some subgroup $G$ of the symmetric group $\Sigma_i$. ${\mathcal I}(i)$ 
is the linear isometries operad (see Section I of~\cite{EKMM}). 
\end{definition}

What makes the class $\F$ of particular interest is the following result, which is listed 
as Theorem 9.5 in~\cite{Maria}; it is a generalization of Theorem VII.6.7 of~\cite{EKMM}. 

\begin{theorem}
For each $M$ in $\F$ let $\Gamma M \to M$ be a cell $A$-module approximation. Then for any finite 
collection $\{M_1, \cdots , M_n \}$ of $\F$, the induced maps 
$$ \Gamma M_1 \wedge \cdots \wedge \Gamma M_n \to M_1 \wedge \cdots \wedge  M_n$$
are weak equivalences of $A$-modules. 
\end{theorem}

In other words, for objects in $\F$ the derived smash product is defined on point set level. 

\begin{definition}
(see Definition 9.6 of~\cite{Maria}) An extended cell is a pair of the form 
($X \wedge B^n _+$, $X \wedge S^{n-1} _+)$, where $n \geq 0$ and 
$A \wedge_S S \wedge_{\mathcal I} {\mathcal I}(i) \ltimes _G K$  for a $G$-spectrum $K$ indexed 
on $U^i$ and which has the homotopy type of a $G$-CW spectrum for some $G<\Sigma_i$. 

An extended cell $A$-module is an $A$-module $M \colim M_i$, with $M_0 =\ast$ and $M_n$ derived 
from $M_{n-1}$ via a pushout of $A$-modules
$$
\xymatrix{
\bigvee_{\alpha} X_\alpha \wedge S^{n_{\alpha}-1} _+ \ar[r] \ar[d]
&
M_{n-1} \ar[d]\\
\bigvee_{\alpha} X_\alpha \wedge B^{n_{\alpha}} _+ \ar[r]
&
M_n.
}
$$
\end{definition}

We intend to show that cell $\a$-algebras are extended cell $A$-modules. For now, we list some 
of the key properties of extended cell modules which make our interest in them evident. 

{\bf 1.} For an extended cell $A$-module $M$ and a subgroup $H$ of the symmetric group 
$\Sigma_n$, the operadic power ${\a}(n) \wedge_{H} M^{\wedge n}$ is in $\F$.

{\bf 2.} For an extended cell $A$-module $M$, a group $H< \Sigma_n$, and an $H$-simplicial 
set $P$, the projection
$$P \wedge_{hH} M^{\wedge n} \to  P \wedge_{H} M^{\wedge n}$$
is a weak equivalence of spectra. 

{\bf 3.} Let $Y \to M$ be a cell $A$-module approximation of the extended cell module $M$. Then 
for all $H< \Sigma_i$ and $H$-simplicial sets $P$,
$$ P \wedge_{H} Y^{\wedge i} \to P \wedge_{H} M^{\wedge i}$$
is a weak equivalence. Moreover, there is a cell $A$-module approximation 
$Z \to P \wedge_{H} M^{\wedge i}$ such that for all $G< \Sigma_j$ and $G$-simplicial sets $Q$,
$$Q \wedge_{G} Z^{\wedge j}  \to 
Q \wedge_{G} (P \wedge_{H} M^{\wedge i})^{\wedge j}$$
is a weak equivalence. 

{\bf 4.} Let $M = P \wedge_H K^{\wedge k}$ and $N = Q \wedge_J L^{\wedge l}$ be $A$-modules
for some extended cell $A$-modules $K$ and $L$, groups $H< \Sigma_k$, $J< \Sigma_l$, and an 
$H$-simplicial set $P$ and $J$-simplicial set $Q$. Then there is a cell $A$-module approximation
$Z \to M \vee N$ such that for all $\Sigma_i$-simplicial sets $T$, 
$$T \wedge_{\Sigma_i} Z^{\wedge i} \to T \wedge_{\Sigma_i} (M \vee N)^{\wedge i}$$
is a weak equivalence.

Properties {\bf 1} and {\bf 2} are consequences of the proof of Theorem 9.8 of~\cite{Maria}, 
Property {\bf 3} follows from the proof of Theorem 9.10 of~\cite{Maria}, and Property {\bf 4}
from Proposition 9.11 once we recall that the 
operad $\a$ arises from an operad in simplicial sets.  

The following (promised) result allows us to take advantage of nice homotopical properties 
of extended cell $A$-modules.

\begin{lemma}
\label{lem:extcel}
Let $C$ be a cell $\a$-algebra. Then it is an extended cell $A$-module. 
\end{lemma}

\begin{proof}
Since $C$ is a cell $\a$-algebra, it can be expressed as $C = colim M_i$, with $M_0=\ast$ 
and $M_i$ obtained from $M_{i-1}$ as a pushout
\begin{equation}
\label{diag:cell}
\xymatrix{
T_{\a} E \ar[r] \ar[d]
&
M_{i-1} \ar[d]\\
T_{\a} CE \ar[r] \ar[r]
&
M_i,
}
\end{equation}
where $E$ is a wedge sphere $A$-modules. 

Consider the simplicial object $\beta_\ast (T_{\a} CE,T_{\a} E,M_{i-1})$ whose 
$p$-simplices are $T_{\a} CE \amalg (T_{\a} E)^{\amalg p} \amalg M_{i-1}$, where
$\amalg$ is the coproduct in the category of $\a$-algebras. The face and degeneracy 
operators are defined using the maps 
$$\mu: T_{\a} CE \amalg T_{\a} E \to T_{\a} \hspace{50pt} and \hspace{50pt}
  \nu: T_{\a} E \amalg M_{i-1} \to M_{i-1},$$
which are induced by the maps in Diagram~\ref{diag:cell}.
See Definition VII.3.5 of~\cite{EKMM} for details. 

A simple comparison of face and degeneracy operators shows that we have an equivalence 
of simplicial objects 
$\beta_\ast (T_{\a} E,T_{\a} E, T_{\a} E)$ and $(T_{\a} E) \otimes I$, where $I$ is the 
standard simplicial $1$-simplex with $p+2$ $p$-simplices. Consequently,

$$\beta_\ast (T_{\a} CE,T_{\a} E,M_{i-1}) \simeq 
T_{\a} CE \amalg_{T_{\a} E} (T_{\a} E \otimes I) \amalg_{T_{\a} E} M_{i-1}.$$
Further, observe that the collapse map induced by $I \to \{pt \}$ 
$$T_{\a} CE \amalg_{T_{\a} E} (T_{\a} E \otimes I) \amalg_{T_{\a} E} M_{i-1} \to 
T_{\a} CE \amalg_{T_{\a} E} M_{i-1} = M_i$$
is an equivalence. This is an immediate consequence of the of a simple homeomorphism 
of based spaces:
$$CX \cup_X (X \wedge I_+) \to CX.$$
The detailed proof of this statement is identical to that of Proposition VII.3.8 of~\cite{EKMM}. 
Thus, 
$$M_i \simeq |\beta_\ast (T_{\a} CE,T_{\a} E,M_{i-1})|.$$ 
To complete the proof, we note that the argument used to prove Lemma VII.7.5 of~\cite{EKMM}
shows that the $q$'th filtration of $\beta_\ast (T_{\a} CE,T_{\a} E,M_{i-1})$ (and by passage 
to colimits any cell $\a$-algebra) is an extended cell $A$-module. 

\end{proof}

Our next objective is to show that $T_{\a} ^k C$ is a sum of objects of the form 
$(P \wedge C^{\wedge n})_H$ where $H$ is a subgroup of the symmetric group $\Sigma_n$ and 
$P$ is an $H$-simplicial set. Observe that if this is the case, then since $C$ is an extended 
cell module by Lemma~\ref{lem:extcel}, Theorem~\ref{theor:tech} is immediate from Properties 
{\bf 1} - {\bf 4}.

As it often happens, the case $k=2$, i.e. $T_{\a} ^2 C$, conveys the essence of the problem, so we 
describe it in detail. 

\begin{eqnarray*}
\lefteqn{
T_{\a} ^2 C = T_{\a} (T_{\a} C) = 
\bigoplus_{i=1} ^{\infty} {\a}(i) \wedge_{\Sigma_i} (T_{\a} C)^{\wedge i}}\\
&\cong& 
\bigoplus _{n=1} ^{\infty} \bigoplus _{l=n} ^{\infty} {\a}(n) \wedge_{\Sigma_n}
[\bigoplus_{i_1 + \cdots + i_n=l} {\a}(i_1) \wedge \cdots \wedge {\a}(i_n) 
\wedge_{\Sigma_{i_1} \times \cdots \times \Sigma_{i_n}} C^{\wedge l}]\\
&\cong&
\bigoplus _{n=1} ^{\infty} \bigoplus _{l=n} ^{\infty} 
\bigoplus_{i_1 + \cdots + i_n=l} 
[{\a}(n) \wedge_{\Sigma_n} {\a}(i_1) \wedge \cdots \wedge {\a}(i_n)
\wedge_{\Sigma_{i_1} \times \cdots \times \Sigma_{i_n}} C^{\wedge l}]\\
&\cong&
\bigoplus _{n=1} ^{\infty} \bigoplus _{l=n} ^{\infty} 
\bigoplus_{i_1 + \cdots + i_n=l}
[{\a}(n) \wedge {\a}(i_1) \wedge \cdots \wedge {\a}(i_n) \wedge 
C^{\wedge l}]_{\Sigma_n \ltimes(\Sigma_{i_1} \times \cdots \times \Sigma_{i_n})},
\end{eqnarray*}
where $\Sigma_n \ltimes(\Sigma_{i_1} \times \cdots \times \Sigma_{i_n})$ is the semi-direct product, 
and is consequently a subgroup of $\Sigma_l$. Thus, each summand is indeed of the form 
$(P \wedge C^{\wedge l})_H$ with $H< \Sigma_l$.

The general case follows by induction. Indeed, $T_{\a} ^{k+1} C = T_{\a} ^{k} (T_{\a} C)$. Assuming 
that $T_{\a} ^{k} C$ is a sum of objects of form $(P \wedge C^{\wedge l})_H$, we need to understand
$(P \wedge (T_{\a} C)^{\wedge l})_H$. Observe that
$$(T_{\a} C)^{\wedge l} = \bigoplus_{(j_1, \cdots , j_l)} {\a}(j_1) \wedge \cdots \wedge {\a}(j_l)
\wedge_{\Sigma_{j_1} \times \cdots \times \Sigma_{j_l}} C^{\wedge (j_1 + \cdots + j_l)}$$ 
Consequently, each summand of $(P \wedge (T_{\a} C)^{\wedge l})_H$ is itself a sum of objects 
of the desired form.

\end{document}